\documentclass{amsart}

\usepackage{amsfonts}
\usepackage{amsmath,amssymb}

\def\Xint#1{\mathchoice
   {\XXint\displaystyle\textstyle{#1}}
   {\XXint\textstyle\scriptstyle{#1}}
   {\XXint\scriptstyle\scriptscriptstyle{#1}}
   {\XXint\scriptscriptstyle\scriptscriptstyle{#1}}
   \!\int}
\def\XXint#1#2#3{{\setbox0=\hbox{$#1{#2#3}{\int}$}
     \vcenter{\hbox{$#2#3$}}\kern-.5\wd0}}

\def\dashint{\Xint-}

\newcommand{\bR}{\mathbb R}
\newcommand{\ip}[1]{\left\langle#1\right\rangle}
\newcommand{\norm}[1]{\left\lVert#1\right\rVert}
\newcommand{\Norm}[1]{\lVert#1\rVert}
\newcommand{\abs}[1]{\left\lvert#1\right\rvert}
\newcommand{\Abs}[1]{\lvert#1\rvert}
\newcommand{\set}[1]{\left\{#1\right\}}

\newcommand{\VMO}{\mathrm{VMO}}
\newcommand{\Hloc}{\mathrm{(H)}_{loc}}

\renewcommand{\vec}[1]{\boldsymbol{#1}}
\DeclareMathOperator*{\osc}{osc}
\DeclareMathOperator{\dist}{dist}

\newtheorem{theorem}{Theorem}[section]
\newtheorem{lemma}[theorem]{Lemma}
\newtheorem{corollary}[theorem]{Corollary}

\theoremstyle{definition}
\newtheorem{definition}[theorem]{Definition}

\theoremstyle{remark}
\newtheorem{remark}[theorem]{Remark}

\numberwithin{equation}{section}

\begin{document}
\title[Green function estimates]{The Green function estimates for strongly elliptic systems of second order}

\author[S. Hofmann]{Steve Hofmann}
\address[S. Hofmann]{Mathematics Department, University of Missouri, Columbia, Missouri 65211, United States of America}
\email{hofmann@math.missouri.edu}

\author[S. Kim]{Seick Kim}
\address[S. Kim]{Centre for Mathematics and its Applications,
The Australian National University, ACT 0200, Australia}
\email{seick.kim@maths.anu.edu.au}




\begin{abstract} We establish existence and pointwise estimates of fundamental solutions
and Green's matrices for divergence form, second order strongly elliptic systems
in a domain $\Omega \subseteq \mathbb{R}^n, n \geq 3,$ under the assumption that solutions of the system
satisfy De Giorgi-Nash type local H\"{o}lder continuity estimates.  In particular, our results apply to perturbations of diagonal systems, and thus especially to complex perturbations of a single real equation.
\end{abstract}
\maketitle

\section{Introduction}\label{sec:I}
In this article, we study Green's functions (or Green's matrices) of second
order, strongly elliptic systems of divergence type in a domain
$\Omega\subset\bR^n$ with $n\ge 3$.  In particular, we treat the Green matrix in the entire space,
usually called the fundamental solution.
We shall prove that if a given elliptic system has
the property that all weak solutions of the system are locally
H\"older continuous, then it has the Green's matrix in $\Omega$.
(For example, if coefficients of the system belong to the space of $\VMO$
introduced by Sarason \cite{Sarason}, then it will enjoy such a property).
For such elliptic systems, we study standard properties of the Green's matrix
including pointwise bounds, $L^p$ and weak $L^p$ estimates for Green's matrix and
its derivatives, etc.

For the scalar case, i.e., a single elliptic equation,
the existence and properties of Green's function was
studied by Littman, Stampacchia, and Weinberger \cite{LSW} and
Gr\"uter and Widman \cite{GW}.
In this article, we follow the approach of Gr\"uter and Widman
in constructing Green's matrix.
The main technical difficulties arise from lack of Harnack type inequalities
and the maximum principle for the systems.
The key observation on which this article is based is that even in
the scalar case, one can get around Moser's Harnack inequality \cite{Moser}
or maximum principle but instead rely solely on De Giorgi-Nash type oscillation
estimates \cite{DeG57} in constructing and studying properties
of Green's functions.
From this point of view, this article provides a unified approach in
studying Green's function for both scalar and systems of equations.
We should point out that there has been some study of Green's matrix for systems with continuous coefficients, notably by Fuchs \cite{Fuchs} and Dolzmann-M\"{u}ller \cite{DM}.
Our existence results and interior estimates of Green's function will include theirs, since as is well known, weak solutions of systems with uniformly continuous (or VMO) coefficients enjoy local H\"older estimates.   On the other hand, we have not attempted to replicate their boundary estimates, which depend in particular on having a $C^1$ boundary.
Our method does not require boundedness of the domain nor regularity of the boundary in constructing Green's matrices, while the methods of Fuchs \cite{Fuchs} and Dolzmann-M\"uller \cite{DM} require both boundedness and regularity of the domain at the very beginning.
We note that a scalar elliptic equation
with complex coefficients can be identified as an elliptic system
with real coefficients satisfying a special structure, and thus our results apply in particular to complex perturbations of
a scalar real equation. 
In the complex coefficients setting, the main results of Section~\ref{sec:E}
in our paper can be also obtained by following the method
of Auscher \cite{Auscher}.  The estimates of the present paper will be applied
to the development of the layer potential method for equations with complex coefficients in
\cite{AAAHK}.

The organization of this paper is as follows.
In Section~\ref{sec:P}, we define the property (H), which is essentially
equivalent to De Giorgi's oscillation estimates in the scalar case,
and introduce a function space $Y^{1,2}_0(\Omega)$
which substitutes $W^{1,2}_0(\Omega)$ in
constructing Green's functions;
they are identical if $\Omega$ is bounded but in general,
$Y^{1,2}_0(\Omega)$ is a larger space and is more suitable for our purpose.
In Section~\ref{sec:E}, we study Green's functions defined in the entire space,
which are usually referred to as the fundamental solutions.
The main result is that for a system whose coefficients are close to those of
a diagonal system, the fundamental solution behaves very much like that
of a single equation.
In Section~\ref{sec:G}, we study Green's matrices in general domains, including
unbounded ones. We also study the boundary behavior of Green's matrices
when the boundary of domain satisfies a measure theoretic
exterior cone condition, called the condition (S).
We prove in particular that if the coefficients of the system are close to those of
a diagonal system, then again the boundary behavior of its Green's function
is much like that of a single equation.
In section~\ref{sec:V}, we discuss the Green's matrices
of the strongly elliptic systems with $\VMO$ coefficients.
By following the same techniques already developed in the
previous two sections, we construct the Green's matrix in general domains
including the entire space.
One subtle difference is that in this $\VMO$ coefficients
case, one should play with a localized version of property (H) since
basically, the regularity of weak solutions of the systems
with $\VMO$ coefficients is inherited from the systems with
constant coefficients when the scale is made small enough.
Therefore, all the estimates for the Green's matrix stated in this
section are only meaningful near a pole.

Finally, we would like to mention that when $n=2$, the method used
in this article breaks down in several places and for that reason
we plan to treat the two dimensional case in a separate paper.

\section{Preliminaries}\label{sec:P}
\subsection{Strongly elliptic systems}
Throughout this article, the summation convention over repeated indices
shall be assumed.
Let $L$ be a second order elliptic operator of divergence type 
acting on vector valued functions
$\vec{u}=(u^1,\ldots,u^N)^T$ defined on $\bR^n$ ($n\ge 3$) in the 
following way:
\begin{equation}
\label{eqP-01}
L\vec{u} = -D_\alpha (\vec{A}^{\alpha\beta}\, D_\beta \vec{u}),
\end{equation}
where $\vec{A}^{\alpha\beta}=\vec{A}^{\alpha\beta}(x)$
($\alpha,\beta=1,\ldots, n$)
are $N$ by $N$ matrices satisfying the strong ellipticity condition, i.e., 
there is a number $\lambda>0$ such that 
\begin{equation}
\label{eqP-02}
A^{\alpha\beta}_{ij}(x)\xi^j_\beta\xi^i_\alpha
\ge \lambda \abs{\vec{\xi}}^2
:=\lambda\sum_{i=1}^N\sum_{\alpha=1}^n\Abs{\xi^i_\alpha}^2,
\quad\forall x\in\bR^n
\end{equation}
We also assume that $A^{\alpha\beta}_{ij}$ are bounded, i.e.,
there is a number $\Lambda>0$ such that 
\begin{equation}
\label{eqP-03}
\sum_{i,j=1}^N\sum_{\alpha,\beta=1}^n
\Abs{A^{\alpha\beta}_{ij}(x)}^2\le \Lambda^2,\quad\forall x\in\bR^n.
\end{equation}
If we write \eqref{eqP-01} component-wise, then we have
\begin{equation}
\label{eqP-04}
(L\vec{u})^i = -D_\alpha (A^{\alpha\beta}_{ij} D_\beta u^j),
\quad \forall i=1,\ldots,N.
\end{equation}
The transpose operator of ${}^t\!L$ of $L$ is defined by
\begin{equation}
\label{eqP-05}
{}^t\!L\vec{u} = -D_\alpha ({}^t\!\vec{A}^{\alpha\beta} D_\beta \vec{u}),
\end{equation}
where ${}^t\!\vec{A}^{\alpha\beta}=(\vec{A}^{\beta\alpha})^T$ 
(i.e., ${}^t\!A^{\alpha\beta}_{ij}=A^{\beta\alpha}_{ji}$).
Note that the coefficients ${}^t\!A^{\alpha\beta}_{ij}$ satisfy
\eqref{eqP-02}, \eqref{eqP-03} with the same constants $\lambda, \Lambda$.

In the sequel, we shall use the notation
$\dashint_{S} f:=\frac{1}{\abs{S}} \int_{S} f$ (assuming $0<|S|<\infty$),
where $S$ is a measurable subset of $\bR^n$ and
$\abs{S}$ denotes the Lebesgue measure of measurable $S$.
\begin{definition}\label{def:P-01}
We say that the operator $L$ satisfies the property (H) if
there exist $\mu_0, H_0>0$ such that all weak solutions $\vec{u}$ of
$L\vec{u}=0$ in $B_R=B_R(x_0)$ satisfy
\begin{equation}
\label{eqP-06}
\int_{B_r}\abs{D \vec{u}}^2 \le H_0 \left(\frac{r}{s}\right)^{n-2+2\mu_0}
\int_{B_s}\abs{D \vec{u}}^2, \quad 0<r<s\le R.
\end{equation}
Similarly, we say that the transpose operator ${}^t\!L$
satisfies the property (H) if corresponding
estimates hold for all weak solutions $\vec{u}$ of
${}^t\!L\vec{u}=0$ in $B_R$.
\end{definition}

\begin{lemma} \label{lem:P-02}
Let $(a^{\alpha\beta}(x))_{\alpha,\beta=1}^n$ be coefficients satisfying
the following conditions:
There are constants $\lambda_0, \Lambda_0>0$ such that for all $x\in\bR^n$
\begin{equation}
\label{eqP-07}
a^{\alpha\beta}(x)\xi_\beta\xi_\alpha\ge
\lambda_0\abs{\xi}^2,\quad\forall\xi\in\bR^n;\quad
\sum\abs{a^{\alpha\beta}(x)}^2\le \Lambda_0^2.
\end{equation}
Then, there exists $\epsilon_0=\epsilon_0(n,\lambda_0,\Lambda_0)$
such that if
\begin{equation}
\label{eqP-08}
\epsilon^2(x):=
\sum\abs{A^{\alpha\beta}_{ij}(x)-a^{\alpha\beta}(x)\delta_{ij}}^2< \epsilon_0^2,
\quad\forall x\in\bR^n,
\end{equation}
then the operator $L$ associated with the coefficients
$A^{\alpha\beta}_{ij}$ satisfies the condition (H)
with $\mu_0=\mu_0(n,\lambda_0,\Lambda_0), H_0=H_0(n,N,\lambda_0,\Lambda_0)>0$.
\end{lemma}
\begin{proof}
See e.g., \cite[Proposition 2.1]{HK}.
\end{proof}

\begin{lemma}\label{lem:P-03}
Suppose that the operator $L$ satisfies the following  
H\"older property for weak solutions:
There are constants $\mu_0, C_0>0$
such that all weak solutions $\vec{u}$ of $L\vec{u}=0$
in $B_{2R}=B_{2R}(x_0)$ satisfy the estimate
\begin{equation}
\label{eqP-09}
[\vec{u}]_{C^{\mu_0}(B_R)}
\le C_0 R^{-\mu_0}\left(\dashint_{B_{2R}}\abs{\vec{u}}^2\right)^{1/2},
\end{equation}
where $[f]_{C^\mu(\Omega)}$ denotes the usual $C^\mu(\Omega)$ semi-norm
of $f$; see \cite{GT} for the definition.
Then, the operator $L$ satisfies the property (H) with
$\mu_0$ and $H_0=H_0(n,N,\lambda,\Lambda,C_0)$.
\end{lemma}
\begin{proof}
We may assume that $r<s/4$; otherwise, \eqref{eqP-06} is trivial.
Denote $\overline{\vec{u}}_r=\dashint_{B_r}\vec{u}$.
We may assume, by replacing $\vec{u}$ by $\vec{u}-\overline{\vec{u}}_{s}$,
if necessary, that $\overline{\vec{u}}_{s}=0$.
From the Caccioppoli inequality, \eqref{eqP-09}, and then
the Poincar\'e inequality, it follows
\begin{equation*}
\begin{split}
\int_{B_r}\abs{D\vec{u}}^2&
\le C r^{-2}\int_{B_{2r}} \abs{\vec{u}-\overline{\vec{u}}_{2r}}^2\le
C r^{-2}\int_{B_{2r}}\dashint_{B_{2r}}
\abs{\vec{u}(x)-\vec{u}(y)}^2\,dy\,dx\\
&\le C r^{-2}[\vec{u}]_{C^{\mu_0}(B_{2r})}^2 (2r)^{2\mu_0}\abs{B_{2r}}
\le C r^{n-2+2\mu_0}[\vec{u}]_{C^{\mu_0}(B_{s/2})}^2\\
& \le C(r/s)^{n-2+2\mu_0}
s^{-2}\int_{B_{s}}\abs{\vec{u}}^2
\le C(r/s)^{n-2+2\mu_0} \int_{B_{s}}\abs{D\vec{u}}^2.
\end{split}
\end{equation*}
The proof is complete.
\end{proof}

\begin{lemma}
\label{lem:P-04}
Assume that the operator $L$ satisfies the property (H).
Then, the operator $L$ satisfies the H\"older property \eqref{eqP-09}.
Moreover, for any $p>0$, there exists
$C_p=C_p(n,N,\lambda,\Lambda,\mu_0,H_0,p)>0$ such that
all weak solutions $\vec{u}$ of $L\vec{u}=0$ in $B_R=B_R(x_0)$ satisfy
\begin{equation}
\label{eqP-11}
\norm{\vec{u}}_{L^\infty(B_r)}\le \frac{C_p}{(R-r)^{n/p}}
\norm{\vec{u}}_{L^p(B_R)},\quad \forall r \in (0,R).
\end{equation}
\end{lemma}
\begin{proof}
From a theorem of Morrey \cite[Thoerem~3.5.2]{Morrey},
the property (H), and the Caccioppoli inequality, it follows that
\begin{equation}
\label{eqP-12}
[\vec{u}]_{C^{\mu_0}(B_R)}^2\le
C R^{2-n-2\mu_0}\norm{D\vec{u}}_{L^2(B_{3R/2})}^2\le
C R^{-n-2\mu_0}\norm{\vec{u}}_{L^2(B_{2R})}^2.
\end{equation}
Then, by a well known averaging argument (see e.g., \cite{HK})
we derive
\begin{equation} \label{eq4.03}
\norm{\vec{u}}_{L^\infty(B_{R/2})}
\leq C \left(\dashint_{B_{2R}}\abs{\vec{u}}^2\right)^{1/2},
\end{equation}
where $C=C(n,N,\lambda,\Lambda,\mu_0,H_0)>0$.
For the proof that \eqref{eq4.03} implies \eqref{eqP-11}, 
we refer to \cite[pp. 80--82]{Gi93}.
\end{proof}

\subsection{Function spaces $Y^{1,2}(\Omega)$ and $Y^{1,2}_0(\Omega)$}
\begin{definition} \label{def:P-05}
For an open set $\Omega\subset\bR^n$ ($n\ge 3$), the
space $Y^{1,2}(\Omega)$ is defined as the family of all weakly differentiable
functions $u\in L^{2^*}(\Omega)$, where $2^*=\frac{2n}{n-2}$,
whose weak derivatives are functions in $L^2(\Omega)$.
The space $Y^{1,2}(\Omega)$ is endowed with the norm
\begin{equation*}
\norm{u}_{Y^{1,2}(\Omega)}:=\norm{u}_{L^{2^*}(\Omega)}+\norm{D u}_{L^2(\Omega)}.
\end{equation*}
We define $Y^{1,2}_0(\Omega)$ as the closure of $C^\infty_c(\Omega)$
in $Y^{1,2}(\Omega)$,
where $C^\infty_c(\Omega)$ is the set of all infinitely differentiable 
functions with compact supports in $\Omega$.
\end{definition}
We note that in the case $\Omega = \mathbb{R}^n$, it is well known that 
$Y^{1,2}(\bR^n)=Y^{1,2}_0(\bR^n)$;
see e.g., \cite[p. 46]{MZ}.
By the Sobolev inequality, it follows that
\begin{equation}
\label{eqP-14}
\norm{u}_{L^{2^*}(\Omega)} \le C(n) \norm{D u}_{L^2(\Omega)},\quad
\forall u\in Y^{1,2}_0(\Omega).
\end{equation}
Therefore, we have $W^{1,2}_0(\Omega)\subset Y^{1,2}_0(\Omega)$
and $W^{1,2}_0(\Omega)=Y^{1,2}_0(\Omega)$ when $\Omega$ has
finite Lebesgue measure.

From \eqref{eqP-14}, it follows that the bilinear form
\begin{equation}
\label{eqP-15}
\ip{\vec{u},\vec{v}}_{\mathbf H}:=\int_{\Omega} D_\alpha u^i D_\alpha v^i
\end{equation}
defines an inner product on $\mathbf H:=Y^{1,2}_0(\Omega)^N$.
Also, it is routine to check that $\mathbf H$ equipped
with the inner product \eqref{eqP-15} is a Hilbert space.
\begin{definition}\label{den:P-06}
We shall denote by $\mathbf H$ the Hilbert space $Y^{1,2}_0(\Omega)^N$
with the inner product \eqref{eqP-15}.
We denote
\begin{equation*}
\norm{\vec{u}}_{\mathbf H}:=\ip{\vec{u},\vec{u}}_{\mathbf H}^{1/2}=\norm{D\vec{u}}_{L^2(\Omega)}.
\end{equation*}
We also define the bilinear form associated to the operator $L$ as
\begin{equation*}
B(\vec{u},\vec{v}):=\int_{\Omega} A^{\alpha\beta}_{ij}
D_\beta u^j D_\alpha v^i.
\end{equation*}
\end{definition}
By the strong ellipticity \eqref{eqP-02}, it follows that
the bilinear form $B$ is coercive; i.e,
\begin{equation}
\label{eqP-18}
B(\vec{u},\vec{u})\ge \lambda \ip{\vec{u},\vec{u}}_{\mathbf H}.
\end{equation}

\section{Fundamental matrix in $\bR^n$}\label{sec:E}
Throughout this section, we assume that the operators
$L$ and ${}^t\!L$ satisfy the property (H).
The main goal of this section is to construct the fundamental matrix
of the the operator $L$ in the entire $\bR^n$, where $n\ge 3$.
Since $Y^{1,2}(\bR^n)=Y^{1,2}_0(\bR^n)$,
we have, as in Definition~\ref{def:P-05},
\begin{equation*}
\norm{u}_{L^{2^*}(\bR^n)} \le C(n) \norm{D u}_{L^{2}(\bR^n)},
\quad \forall u\in Y^{1,2}(\bR^n).
\end{equation*}
We note that $W^{1,2}(\bR^n)\subset Y^{1,2}(\bR^n)\subset W^{1,2}_{loc}(\bR^n)$.
Unless otherwise stated, we employ the letter $C$ to denote a constant
depending on $n$, $N$, $\lambda$, $\Lambda$, $\mu_0$, $H_0$, and sometimes 
on an exponent $p$ characterizing Lebesgue classes.
It should be understood that $C$ may vary from line to line.

\subsection{Averaged fundamental matrix}\label{sec:E-01}
Our approach here is based on that in \cite{GW}.
Let $y\in\bR^n$ and $1\le k\le N$ be fixed.
For $\rho>0$, consider the linear functional 
$\vec{u}\mapsto \dashint_{B_\rho(y)} u^k$.
Since
\begin{equation}
\label{eqE-02}
\abs{\,\dashint_{B_\rho(y)}u^k\,}\le
C\rho^{(2-n)/2}\norm{\vec{u}}_{L^{2^*}(\bR^n)}\le
C\rho^{(2-n)/2}\norm{\vec{u}}_{\mathbf H},
\end{equation}
Lax-Milgram lemma implies that there exists
a unique $\vec{v}_\rho=\vec{v}_{\rho;y,k}\in {\mathbf H}$ such that
\begin{equation}
\label{eqE-03}
\int_{\bR^n}A^{\alpha\beta}_{ij}D_\beta v^j_{\rho}\,D_\alpha u^i=
\dashint_{B_\rho(y)}u^k,\quad\forall\vec{u}\in {\mathbf H}.
\end{equation}
Note that \eqref{eqP-18}, \eqref{eqE-03}, and \eqref{eqE-02} imply that
\begin{equation*}
\lambda \norm{\vec{v}_\rho}_{\mathbf H}^2\le
B(\vec{v}_\rho,\vec{v}_\rho)\le C\rho^{(2-n)/2} \norm{\vec{v}_\rho}_{\mathbf H},
\end{equation*}
and thus we have
\begin{equation}
\label{eqE-05}
\norm{D\vec{v}_\rho}_{L^2(\bR^n)}=\norm{\vec{v}_\rho}_{\mathbf H}\le C\rho^{(2-n)/2}.
\end{equation}
We define the ``averaged fundamental matrix''
$\vec{\Gamma}^\rho(\,\cdot\,,y)=(\Gamma_{jk}^\rho(\,\cdot\,,y))_{j,k=1}^N$ by
\begin{equation}
\label{eqE-06}
\Gamma^\rho_{jk}(\,\cdot\,,y)=v^j_{\rho}=v^j_{\rho;y,k}.
\end{equation}
Note that we have
\begin{equation}
\label{eqE-07}
\int_{\bR^n}A^{\alpha\beta}_{ij} D_\beta \Gamma^\rho_{jk}(\,\cdot\,,y)
D_\alpha u^i = \dashint_{B_\rho(y)}u^k,
\quad\forall\vec{u}\in {\mathbf H},
\end{equation}
and equivalently ($\alpha\leftrightarrow\beta$ , $i\leftrightarrow j$).
\begin{equation}
\label{eqE-08}
\int_{\bR^n}{}^t\!A^{\alpha\beta}_{ij} D_\beta u^j
D_\alpha\Gamma^\rho_{ik}(\,\cdot\,,y)= \dashint_{B_\rho(y)}u^k,
\quad\forall\vec{u}\in {\mathbf H}.
\end{equation}

In the sequel, we shall denote by $L^\infty_c(\Omega)$
the family of all $L^\infty$ functions with compact supports in $\Omega$.
For a given $\vec{f}\in L^\infty_c(\bR^n)^N$
consider a linear functional
\begin{equation}
\label{eqE-09}
\vec{w}\mapsto \int_{\bR^n}\vec{f}\cdot\vec{w},
\end{equation}
which is bounded on ${\mathbf H}$ since
\begin{equation}
\label{eqE-10}
\abs{\,\int_{\bR^n} \vec{f}\cdot\vec{w}\,}\le
\norm{\vec{f}}_{L^{\frac{2n}{n+2}}(\bR^n)}
\norm{\vec{w}}_{L^{\frac{2n}{n-2}}(\bR^n)}
\le C \norm{\vec{f}}_{L^{\frac{2n}{n+2}}(\bR^n)} \norm{\vec w}_{\mathbf H}.
\end{equation}
Therefore, by Lax-Milgram lemma, there exists $\vec{u}\in \mathbf H$ such that
\begin{equation}
\label{eqE-11}
\int_{\bR^n}{}^t\!A^{\alpha\beta}_{ij} D_\beta u^j
D_\alpha w^i= \int_{\bR^n} f^i w^i,
\quad\forall\vec{w}\in \mathbf H.
\end{equation}
In particular, if we set $\vec{w}=\vec{v}_\rho$ in \eqref{eqE-11},
then by \eqref{eqE-08}, we have
\begin{equation}
\label{eqE-12}
\int_{\bR^n} \Gamma^\rho_{ik}(\,\cdot\,,y) f^i=
\dashint_{B_\rho(y)}u^k.
\end{equation}
Moreover, by setting $\vec{w}=\vec{u}$ in \eqref{eqE-11}, it follows
from \eqref{eqE-10} that
\begin{equation}
\label{eqE-13}
\norm{D\vec{u}}_{L^2(\bR^n)}\le C
\norm{\vec{f}}_{L^{2n/(n+2)}(\bR^n)}.
\end{equation}

\subsection{$L^\infty$ estimates for averaged fundamental matrix}
\label{sec:E-02}
Let $\vec{u}\in \mathbf H$ be given as in \eqref{eqE-11}. 
We will obtain local $L^\infty$ estimates for $\vec{u}$ in $B_R(x_0)$, where
$x_0\in\bR^n$ and $R>0$ are fixed but arbitrary.

Fix $x\in B_R(x_0)$ and $0<s\le R$.
We decompose $\vec{u}$ as $\vec{u}=\vec{u}_1+\vec{u}_2$, where
$\vec{u}_1\in W^{1,2}(B_s(x))^N$ is the weak solution of
${}^t\!L\vec{u}_1=0$ in $B_s(x)$ satisfying
$\vec{u}_1=\vec{u}$ on $\partial B_s(x)$; i.e.,
$\vec{u}_1-\vec{u}\in W^{1,2}_0(B_s(x))$.
Then, for $0<r<s$, we have
\begin{equation*}
\begin{split}
\int_{B_r(x)} \abs{D\vec{u}}^2
&\le 2\int_{B_r(x)}\abs{D\vec{u}_1}^2+2\int_{B_r(x)} \abs{D\vec{u}_2}^2\\
&\le C\left(\frac{r}{s}\right)^{n-2+2\mu_0}
\int_{B_s(x)} \abs{D\vec{u}_1}^2+ 2\int_{B_s(x)}\abs{D\vec{u}_2}^2\\
&\le C\left(\frac{r}{s}\right)^{n-2+2\mu_0}
\int_{B_s(x)} \abs{D\vec{u}}^2+C\int_{B_s(x)} \abs{D\vec{u}_2}^2.
\end{split}
\end{equation*}
Since $\vec{u}_2\in W^{1,2}_0(B_s(x))^N$ is a weak solution of
${}^t\!L\vec{u}_2=\vec{f}$ in $B_s(x)$, we have
\begin{equation*}
\int_{B_s(x)} \abs{D\vec{u}_2}^2 \le
C \norm{\vec{f}}_{L^{2n/(n+2)}(B_s(x))}^2.
\end{equation*}
For given $p>n/2$, choose $p_0\in (n/2,p)$ such that $\mu_1:=2-n/p_0<\mu_0$.
Then
\begin{equation}
\label{eqE-16}
\norm{\vec{f}}_{L^{\frac{2n}{n+2}}(B_s(x))}^2\le
\norm{\vec{f}}_{L^{p_0}(B_s(x))}^2 \abs{B_s}^{1+2/n-2/{p_0}}
\le C \norm{\vec{f}}_{L^{p_0}(\bR^n)}^2 s^{n-2+2\mu_1}.
\end{equation}
Therefore, after combining the above inequalities, we have
for all $r<s\le R$
\begin{equation*}
\int_{B_r(x)} \abs{D\vec{u}}^2
\le C\left(\frac{r}{s}\right)^{n-2+2\mu_0}
\int_{B_s(x)} \abs{D\vec{u}}^2+
C s^{n-2+2\mu_1}\norm{\vec{f}}_{L^{p_0}(\bR^n)}^2.
\end{equation*}
By a well known iteration argument
(see e.g., \cite[Lemma~2.1, p. 86]{Gi83}), we have
\begin{equation}
\label{eqE-18}
\begin{split}
\int_{B_r(x)}\abs{D\vec{u}}^2
&\le C\left(\frac{r}{R}\right)^{n-2+2\mu_1}\int_{B_R(x)}\abs{D\vec{u}}^2
+C r^{n-2+2\mu_1}\norm{\vec{f}}_{L^{p_0}(\bR^n)}^2\\
&\le C\left(\frac{r}{R}\right)^{n-2+2\mu_1}\int_{\bR^n}\abs{D\vec{u}}^2
+C r^{n-2+2\mu_1}\norm{\vec{f}}_{L^{p_0}(\bR^n)}^2,
\end{split}
\end{equation}
for all $0<r<R$ and $x\in B_R(x_0)$.
From \eqref{eqE-18} it follows (see, e.g. \cite{HK})
\begin{equation}
\label{eqE-19}
[\vec{u}]^2_{C^{\mu_1}(B_R(x_0))} \le C
\left(R^{-(n-2+2\mu_1)}\norm{D\vec{u}}_{L^2(\bR^n)}^2+
\norm{\vec{f}}_{L^{p_0}(\bR^n)}^2\right).
\end{equation}
Note that since $\vec{u}\in \mathbf H$, we have
\begin{equation*}
\norm{\vec{u}}_{L^2(B_R(x_0))}^2 \le \norm{\vec{u}}_{L^{2^*}(B_R(x_0))}^2
\abs{B_R}^{2/n} \le C R^2\norm{D\vec{u}}_{L^2(\bR^n)}^2.
\end{equation*}
Consequently, we have
\begin{equation*}
\begin{split}
\norm{\vec{u}}_{L^\infty(B_{R/2}(x_0))}^2
&\le C R^{2\mu_1} [\vec{u}]_{C^{\mu_1}(B_R(x_0))}^2
+C R^{-n}\norm{\vec{u}}_{L^2(B_{R}(x_0))}^2 \\
&\le C \left(R^{2-n}\norm{D\vec{u}}_{L^2(\bR^n)}^2+
R^{2\mu_1}\norm{\vec{f}}_{L^{p_0}(\bR^n)}^2\right)
+C R^{2-n}\norm{D\vec{u}}_{L^2(\bR^n)}^2 \\
&\le C R^{2-n}\norm{\vec{f}}_{L^{2n/(n+2)}(\bR^n)}^2+
C R^{2\mu_1}\norm{\vec{f}}_{L^{p_0}(\bR^n)}^2,
\end{split}
\end{equation*}
where we used the inequality \eqref{eqE-13} in the last step.

Therefore, if $\vec{f}$ is supported in $B_R(x_0)$, then \eqref{eqE-16}
yields (recall $\mu_1=2-n/p_0$)
\begin{equation}
\label{eqE-21}
\norm{\vec{u}}_{L^\infty(B_{R/2}(x_0))}
\le C R^{2-n/p_0} \norm{\vec{f}}_{L^{p_0}(B_R(x_0))}
\le C R^{2-n/p} \norm{\vec{f}}_{L^p(B_R(x_0))}.
\end{equation}
Now, \eqref{eqE-12} implies that for $\rho<R/2$,
we have, by setting $x_0=y$ in \eqref{eqE-21},
\begin{equation}
\label{eqE-22}
\abs{\int_{B_R(y)} \Gamma^\rho_{ik}(\,\cdot\,,y) f^i\,}
\le \dashint_{B_\rho(y)}\abs{\vec{u}} \le
C R^{2-n/p} \norm{\vec{f}}_{L^p(B_R(y))},\quad \forall p>n/2
\end{equation}
provided that $\vec{f}$ is supported in $B_R(y)$.
Therefore, by duality, we see that
\begin{equation}
\label{eqE-23}
\norm{\vec{v}_\rho}_{L^q(B_R(y))}\le C R^{2-n+n/q},
\quad\forall q\in [1,\tfrac{n}{n-2}), \quad \forall \rho \in(0,R/2),
\end{equation}
where $\vec{v}_\rho=\vec{v}_{\rho;y,k}$ is as in \eqref{eqE-06}.

Fix $x\neq y$ and let $r:=\frac{2}{3}\abs{x-y}$. If $\rho<r/2$, then
since $\vec{v}_\rho\in W^{1,2}(B_r(x))^N$ and satisfies $L\vec{v}_\rho=0$
weakly in $B_r(x)$, it follows from Lemma~\ref{lem:P-04} that
\begin{equation}
\label{eqE-24}
\abs{\vec{v}_\rho(x)}
\le C r^{-n} \norm{\vec{v}_\rho}_{L^1(B_r(x))}
\le C r^{-n} \norm{\vec{v}_\rho}_{L^1(B_{3r}(y))}\le C r^{2-n}.
\end{equation}
Since $\rho, y, k$ are arbitrary, we have obtained the following estimates.
\begin{equation}
\label{eqE-25}
\abs{\vec{\Gamma}^\rho(x,y)}\le C \abs{x-y}^{2-n},\quad
\forall \rho<\abs{x-y}/3.
\end{equation}

\subsection{Uniform weak-$L^{\frac{n}{n-2}}$ estimates for $\vec{\Gamma}^\rho(\,\cdot\,,y)$}\label{sec:E-03}
We claim that the following estimate holds:
\begin{equation}
\label{eqE-26}
\int_{\bR^n\setminus B_R(y)}
\abs{\vec{\Gamma}^\rho(\,\cdot\,,y)}^{\frac{2n}{n-2}} \le C R^{-n},
\quad \forall R>0,\quad \forall \rho>0.
\end{equation}
If $R> 3\rho$, then by \eqref{eqE-25} we have
\begin{equation*}
\int_{\bR^n\setminus B_R(y)}
\abs{\vec{\Gamma}^\rho(x,y)}^{\frac{2n}{n-2}}\,dx \le C
\int_{\bR^n\setminus B_R(y)} \abs{x-y}^{-2n}\,dx \le C R^{-n}.
\end{equation*}
Next, we consider the case $R\le 3\rho$. Let $\vec{v}_\rho^T$ be
the $k$-th column of the averaged fundamental
matrix $\vec{\Gamma}^\rho(\,\cdot\,,y)$ as in \eqref{eqE-06}.
From \eqref{eqE-05}, we see that
\begin{equation*}
\norm{\vec{v}_\rho}_{L^{2^*}(\bR^n\setminus B_R(y))}
\le \norm{\vec{v}_\rho}_{L^{2^*}(\bR^n)}
\le \norm{D\vec{v}_\rho}_{L^2(\bR^n)} \le C\rho^{(2-n)/2}.
\end{equation*}
and thus \eqref{eqE-26} also follows in the case when $R\le 3\rho$.

Now, let $A_t=\set{x\in\bR^n:\abs{\vec{\Gamma}^\rho(x,y)}>t}$ and
choose $R=t^{-1/(n-2)}$.
Then,
\begin{equation*}
\abs{A_t\setminus B_R(y)}\le t^{-\frac{2n}{n-2}}
\int_{A_t\setminus B_R(y)}
\abs{\vec{\Gamma}^\rho(\,\cdot\,,y)}^{\frac{2n}{n-2}}
\le Ct^{-\frac{2n}{n-2}}\,t^{\frac{n}{n-2}}
=Ct^{-\frac{n}{n-2}}.
\end{equation*}
Obviously, $\abs{A_t\cap B_R(y)}\le C R^n =Ct^{-\frac{n}{n-2}}$.
Therefore, we obtained that for all $t>0$, we have 
\begin{equation}
\label{eqE-30}
\abs{\set{x\in\bR^n:\abs{\vec{\Gamma}^\rho(x,y)}>t}}\le C t^{-\frac{n}{n-2}},
\quad \forall \rho>0.
\end{equation}

\subsection{Uniform weak-$L^{\frac{n}{n-1}}$ estimates for $D\vec{\Gamma}^\rho(\,\cdot\,y)$}
Let $\vec{v}_\rho$ be as before.
Fix a cut-off function $\eta\in C^\infty(\bR^n)$ such that
$\eta\equiv 0$ on $B_{R/2}(y)$,  $\eta\equiv 1$ outside $B_R(y)$, and
$\abs{D\eta}\le C/R$.
If we set $\vec{u}:=\eta^2\vec{v}_\rho$, then by \eqref{eqE-03}
\begin{equation*}
0=\int_{\bR^n}\eta^2A^{\alpha\beta}_{ij}D_\beta v^j_{\rho}D_\alpha v_\rho^i
+\int_{\bR^n}2\eta A^{\alpha\beta}_{ij}D_\beta v^j_{\rho}
v^i_\rho D_\alpha\eta,
\end{equation*}
which together with \eqref{eqE-25} implies that
if $R>6\rho$, then
\begin{equation*}
\int_{\bR^n\setminus B_R(y)}\abs{D \vec{v}_{\rho}}^2
\le C R^{-2}\int_{B_R(y)\setminus B_{R/2}(y)}\abs{\vec{v}_{\rho}}^2
\le C R^{2-n}.
\end{equation*}
On the other hand, if $R\le 6\rho$, then \eqref{eqE-05} again implies
\begin{equation*}
\int_{\bR^n\setminus B_R(y)}\abs{D \vec{v}_{\rho}}^2\le
\int_{\bR^n}\abs{D \vec{v}_{\rho}}^2 \le C \rho^{2-n}
\le C R^{2-n}.
\end{equation*}
Therefore, we have
\begin{equation}
\label{eqE-34}
\int_{\bR^n\setminus B_R(y)}\abs{D \vec{\Gamma}^{\rho}(\,\cdot\,,y)}^2
\le C R^{2-n},\quad\forall R>0,\quad\forall \rho>0.
\end{equation}
Next, let $A_t=\set{x\in\bR^n:\abs{D_x\vec{\Gamma}^\rho(x,y)}>t}$ and
choose $R=t^{-1/(n-1)}$.
Then
\begin{equation*}
\abs{A_t\setminus B_R(y)}\le t^{-2}
\int_{A_t\setminus B_R(y)}\abs{D\vec{\Gamma}^\rho(\,\cdot\,,y)}^2 \le
Ct^{-\frac{n}{n-1}}
\end{equation*}
and $\abs{A_t\cap B_R(y)}\le C R^n =Ct^{-\frac{n}{n-1}}$.
We have thus find that for all $t>0$, we have 
\begin{equation}
\label{eqE-36}
\abs{\set{x\in\bR^n:\abs{D_x \vec{\Gamma}^\rho(x,y)}>t}}\le
C t^{-\frac{n}{n-1}},\quad \forall \rho>0.
\end{equation}

\subsection{Construction of the fundamental matrix}\label{sec:2-5}
First, we claim 
\begin{equation}
\label{eqE-37}
\norm{D \vec{\Gamma}^\rho(\,\cdot\,,y)}_{L^p(B_R(y))}\le C_p R^{1-n+n/p},
\quad\forall \rho>0,\quad \forall p\in (0,\tfrac{n}{n-1}).
\end{equation}
Let $\vec{v}_\rho$ be as before.
Note that
\begin{equation*}
\begin{split}
\int_{B_R(y)} \abs{D\vec{v}_\rho}^p&=
\int_{B_R(y)\cap{\set{\abs{D\vec{v}^\rho}\le \tau}}}\abs{D\vec{v}_\rho}^p+
\int_{B_R(y)\cap{\set{\abs{D\vec{v}^\rho}> \tau}}}\abs{D\vec{v}_\rho}^p\\
&\le  \tau^p\abs{B_R}+
\int_{\set{\abs{D\vec{v}^\rho}> \tau}}\abs{D\vec{v}_\rho}^p.
\end{split}
\end{equation*}
By using \eqref{eqE-36}, we estimate
\begin{equation*}
\begin{split}
\int_{\set{\abs{D\vec{v}^\rho}> \tau}}\abs{D\vec{v}_\rho}^p
&=\int_0^\infty p t^{p-1}\abs{\set{\abs{D\vec{v}^\rho}> \max(t,\tau)}}\,dt\\
&\le C \tau^{-\frac{n}{n-1}}\int_0^\tau p t^{p-1}\,dt
+ C \int_\tau^\infty p t^{p-1-n/(n-1)}\,dt\\
&=C\left(1-p/(p-\tfrac{n}{n-1})\right) \tau^{p-n/(n-1)}.
\end{split}
\end{equation*}
By optimizing over $\tau$, we get
\begin{equation}
\label{eqE-40}
\int_{B_R(y)} \abs{D\vec{v}_\rho}^p\le
C R^{(1-n)p+n},
\end{equation}
from which \eqref{eqE-37} follows.

If we utilize \eqref{eqE-30} instead of \eqref{eqE-36},
we obtain a similar estimates for $\vec{\Gamma}^\rho(\,\cdot\,,y)$
\begin{equation}
\label{eqE-41}
\norm{\vec{\Gamma}^\rho(\,\cdot\,,y)}_{L^p(B_R(y))}\le C_p R^{2-n+n/p},
\quad\forall \rho>0,\quad \forall p\in (0,\tfrac{n}{n-2}).
\end{equation}

Let us fix $q\in (1,\frac{n}{n-1})$.
We have seen that for all $R>0$, there exists some $C(R)<\infty$ such that
\begin{equation*}
\norm{\vec{\Gamma}^\rho(\,\cdot\,,y)}_{W^{1,q}(B_R(y))} \le C(R),
\quad\forall\rho>0.
\end{equation*}
Therefore, by a diagonalization process, we obtain
a sequence $\set{\rho_\mu}_{\mu=1}^\infty$ and
$\vec{\Gamma}(\,\cdot\,,y)$ in $W^{1,q}_{loc}(\bR^n)^{N\times N}$ such that
$\lim_{\mu\to\infty}\rho_\mu=0$ and that
\begin{equation}
\label{eqE-43}
\vec{\Gamma}^{\rho_\mu}(\,\cdot\,,y) \rightharpoonup\vec{\Gamma}(\,\cdot\,,y)
\text{ in } W^{1,q}(B_R(y))^{N\times N}, \quad \forall R>0,
\end{equation}
where we recall that $\rightharpoonup$ denotes weak convergence.
Then, for any $\vec{\phi}\in C^\infty_c(\bR^n)^N$, it follows from
\eqref{eqE-07}
\begin{equation}
\label{eqE-44}
\begin{split}
\int_{\bR^n}A^{\alpha\beta}_{ij} D_\beta \Gamma_{jk}(\,\cdot\,,y)
D_\alpha \phi^i &=
\lim_{\mu\to\infty}\int_{\bR^n}A^{\alpha\beta}_{ij}
D_\beta \Gamma^{\rho_\mu}_{jk}(\,\cdot\,,y) D_\alpha \phi^i \\
&=\lim_{\mu\to\infty}\dashint_{B_{\rho_\mu}(y)}\phi^k=
\phi^k(y).
\end{split}
\end{equation}
Let $\vec{v}_\rho^T$ be the $k$-th column of $\vec{\Gamma}^\rho(\,\cdot\,,y)$
as before, and let $\vec{v}^T$ be the corresponding $k$-th column
of $\vec{\Gamma}(\,\cdot\,,y)$.
Then, for any $\vec{g}\in L^\infty_{c}(B_R(y))^N$, \eqref{eqE-41}
yields 
\begin{equation}
\label{eqE-45}
\abs{\int_{\bR^n}\vec{v}\cdot\vec{g}\,}
=\lim_{\mu\to\infty}\abs{\int_{\bR^n}\vec{v}_{\rho_\mu}\cdot\vec{g}\,}
\le C_pR^{2-n+n/p}\norm{\vec{g}}_{L^{p'}(B_R(y))},
\end{equation}
where $p'$ denotes the conjugate exponent of $p\in [1,\tfrac{n}{n-2})$.
Therefore, we obtain
\begin{equation}
\label{eqE-46}
\norm{\vec{\Gamma}(\,\cdot\,,y)}_{L^p(B_R(y))}\le C_p\,R^{2-n+n/p},
\quad \forall p\in [1,\tfrac{n}{n-2}).
\end{equation}
By a similar reasoning, we also have by \eqref{eqE-37}
\begin{equation}
\label{eqE-47}
\norm{D\vec{\Gamma}(\,\cdot\,,y)}_{L^p(B_R(y))}\le C_p\,R^{1-n+n/p},
\quad \forall p\in [1,\tfrac{n}{n-1}).
\end{equation}
Also, with the aid of \eqref{eqE-26} and \eqref{eqE-34}, we obtain
\begin{gather}
\label{eqE-48}
\int_{\bR^n\setminus B_R(y)}\abs{\vec{\Gamma}(\,\cdot\,,y)}^{2^*}
\le C R^{-n},\\
\label{eqE-49}
\int_{\bR^n\setminus B_R(y)}\abs{D \vec{\Gamma}(\,\cdot\,,y)}^2
\le C R^{2-n}.
\end{gather}
In particular, \eqref{eqE-48}, \eqref{eqE-49} imply that 
\begin{equation}
\label{eqE-50}
\norm{\vec{\Gamma}(\,\cdot\,,y)}_{Y^{1,2}(\bR^n\setminus B_r(y))}
\le C r^{1-n/2}, \quad \forall r>0.
\end{equation}
Moreover, arguing as before, we see that
the estimates \eqref{eqE-48} and \eqref{eqE-49} imply
\begin{gather}
\label{eqE-51}
\abs{\set{x\in\bR^n:\abs{\vec{\Gamma}(x,y)}>t}}\le C t^{-\frac{n}{n-2}},
\quad\forall t>0\\
\label{eqE-52}
\abs{\set{x\in\bR^n:\abs{D_x \vec{\Gamma}(x,y)}>t}}\le
C t^{-\frac{n}{n-1}} \quad\forall t>0.
\end{gather}

Next, we turn to pointwise bounds for $\vec{\Gamma}(\,\cdot\,,y)$.
Let $\vec{v}^T$ be the $k$-th column of $\vec{\Gamma}(\,\cdot\,,y)$.
For each $x\neq y$, denote $r=\tfrac{2}{3}\abs{x-y}$.
Then, it follows from \eqref{eqE-50} and \eqref{eqE-44}
that $\vec{v}$ is a weak solution of $L\vec{v}=0$
in $B_r(x)$.
Therefore, by Lemma~\ref{lem:P-04} and \eqref{eqE-46} we find 
\begin{equation}
\label{eqE-53}
\abs{\vec{v}(x)}\le C r^{-n}\norm{\vec{v}}_{L^1(B_r(x))}
\le C r^{-n}\norm{\vec{v}}_{L^1(B_{3r}(y))}\le C r^{2-n},
\end{equation}
from which it follows
\begin{equation}
\label{eqE-54}
\abs{\vec{\Gamma}(x,y)}\le C\abs{x-y}^{2-n},\quad \forall x\neq y.
\end{equation}

\subsection{Continuity of the fundamental matrix}\label{sec:E-06}
From the property (H),
it follows that $\vec{\Gamma}(\,\cdot\,,y)$ is H\"older continuous in
$\bR^n\setminus\set{y}$.
In fact, \eqref{eqP-12} together with \eqref{eqE-44} and \eqref{eqE-49}
implies
\begin{equation}
\label{eqE-55}
\abs{\vec{\Gamma}(x,y)-\vec{\Gamma}(z,y)}\le
C \abs{x-z}^{\mu_0} \abs{x-y}^{2-n-\mu_0}
\quad\text{if }\abs{x-z}< \abs{x-y}/2.
\end{equation}
Moreover, by the same reasoning, it follows from \eqref{eqP-12} and
\eqref{eqE-34}
that for any given compact set $K\Subset \bR^n\setminus\set{y}$,
the sequence $\set{\vec{\Gamma}^{\rho_\mu}(\,\cdot\,,y)}_{\mu=1}^\infty$
is equicontinuous on $K$.
Also, by Lemma~\ref{lem:P-04} and \eqref{eqE-26}, we find that
there are $C_K<\infty$ and $\rho_K>0$ such that
\begin{equation}
\label{eqE-55b}
\norm{\vec{\Gamma}^\rho(\,\cdot\,,y)}_{L^\infty(K)} \leq C_K
\quad \forall \rho<\rho_K\quad
\text{for any compact } K\Subset \bR^n\setminus\set{y}.
\end{equation}
Therefore, we may assume, by passing if necessary to a subsequence, that
\begin{equation}
\label{eqE-56}
\vec{\Gamma}^{\rho_\mu}(\,\cdot\,,y) \rightarrow\vec{\Gamma}(\,\cdot\,,y)
\text{ uniformly on $K$, for any compact }
K\Subset \bR^n\setminus\set{y}.
\end{equation}

We will now show that $\vec{\Gamma}(x,\,\cdot\,)$ is also H\"older continuous
in $\bR^n\setminus\set{x}$.
Denote by ${}^t\!\vec{\Gamma}^{\sigma}(\,\cdot\,,x)$ the averaged
fundamental matrix associated to ${}^t\!L$, the transpose of $L$.
Since each column of 
$\vec{\Gamma}^{\rho}(\,\cdot\,,y)$ and
${}^t\!\vec{\Gamma}^{\sigma}(\,\cdot\,,x)$ belongs to $\mathbf H$,
we have by \eqref{eqE-07},
\begin{equation}
\label{eqE-59}
\begin{split}
\dashint_{B_{\rho}(y)}{}^t\!\Gamma_{kl}^{\sigma}(\,\cdot\,,x)
&=\int_{\bR^n}A^{\alpha\beta}_{ij} D_\beta \Gamma^{\rho}_{jk}(\,\cdot\,,y)
D_\alpha {}^t\!\Gamma^{\sigma}_{il}(\,\cdot\,,x)\\
&=\int_{\bR^n}{}^t\!A^{\beta\alpha}_{ji}
D_\alpha {}^t\!\Gamma^{\sigma}_{il}(\,\cdot\,,x)
D_\beta \Gamma^{\rho}_{jk}(\,\cdot\,,y)
=\dashint_{B_{\sigma}(x)}\Gamma_{lk}^{\rho}(\,\cdot\,,y).
\end{split}
\end{equation}
By the same argument as appears in Sec.~\ref{sec:2-5},
we obtain a sequence $\set{\sigma_\nu}_{\nu=1}^\infty$ tending to $0$ such that
${}^t\!\vec{\Gamma}^{\sigma_\nu}(\,\cdot\,,x)$ converges to
${}^t\!\vec{\Gamma}(\,\cdot\,,x)$ uniformly on
any compact subset of $\bR^n\setminus\set{x}$, where 
${}^t\!\vec{\Gamma}(\,\cdot\,,x)$ is a fundamental matrix for ${}^t\!L$
satisfying all properties stated in Sec.~\ref{sec:2-5}.
By \eqref{eqE-59}, we find that
\begin{equation*}
g^{kl}_{\mu\nu}:=
\dashint_{B_{\rho_\mu}(y)}{}^t\!\Gamma_{kl}^{\sigma_\nu}(\,\cdot\,,x)
=\dashint_{B_{\sigma_\nu}(x)}\Gamma_{lk}^{\rho_\mu}(\,\cdot\,,y).
\end{equation*}
From the continuity of $\Gamma_{lk}^{\rho_\mu}(\,\cdot\,,y)$, it follows
that for $x, y\in \bR^n$ with $x\neq y$, we have
\begin{equation*}
\lim_{\nu\to\infty} g^{kl}_{\mu\nu}=
\lim_{\nu\to\infty}
\dashint_{B_{\sigma_\nu}(x)}\Gamma_{lk}^{\rho_\mu}(\,\cdot\,,y)=
\Gamma_{lk}^{\rho_\mu}(x,y)
\end{equation*}
and thus by \eqref{eqE-56} we obtain
\begin{equation*}
\lim_{\mu\to\infty} \lim_{\nu\to\infty} g^{kl}_{\mu\nu}=
\lim_{\mu\to\infty} \Gamma_{lk}^{\rho_\mu}(x,y)=\Gamma_{lk}(x,y). 
\end{equation*}
On the other hand, \eqref{eqE-43} yields
\begin{equation*}
\lim_{\nu\to\infty} g^{kl}_{\mu\nu}=
\lim_{\nu\to\infty}
\dashint_{B_{\rho_\mu}(y)}{}^t\!\Gamma_{kl}^{\sigma_\nu}(\,\cdot\,,x)=
\dashint_{B_{\rho_\mu}(y)}{}^t\!\Gamma_{kl}(\,\cdot\,,x)
\end{equation*}
and thus it follows from the continuity of
${}^t\!\Gamma_{kl}(\,\cdot\,,x)$ that
\begin{equation*}
\lim_{\mu\to\infty} \lim_{\nu\to\infty} g^{kl}_{\mu\nu}=
\lim_{\mu\to\infty} \dashint_{B_{\rho_\mu}(y)}{}^t\!\Gamma_{kl}(\,\cdot\,,x)
={}^t\!\Gamma_{kl}(y,x).
\end{equation*}
We have thus shown that 
\begin{equation*}
\Gamma_{lk}(x,y)={}^t\!\Gamma_{kl}(y,x),
\quad \forall k,l=1,\ldots,N,\quad \forall x\neq y,
\end{equation*}
which is equivalent to say
\begin{equation}
\label{eqE-66}
\vec{\Gamma}(x,y)={}^t\!\vec{\Gamma}(y,x)^T,\quad \forall x\neq y.
\end{equation}
Therefore, we have proved the claim that $\vec{\Gamma}(x,\,\cdot\,)$
is H\"older continuous in $\bR^n\setminus\set{x}$.

So far, we have seen that there is a sequence
$\set{\rho_\mu}_{\mu=1}^\infty$ tending to $0$
such that $\vec{\Gamma}^{\rho_\mu}(\,\cdot\,,y) \to \vec{\Gamma}(\,\cdot\,,y)$
in $\bR^n\setminus\set{y}$.
However, by \eqref{eqE-59}, we obtain
\begin{equation}
\label{eqE-67}
\begin{split}
\Gamma_{lk}^{\rho}(x,y)&=
\lim_{\nu\to\infty}\dashint_{B_{\sigma_\nu}(x)}\Gamma_{lk}^{\rho}(\,\cdot\,,y)
=\lim_{\nu\to\infty}
\dashint_{B_{\rho}(y)}{}^t\!\Gamma_{kl}^{\sigma_\nu}(\,\cdot\,,x)\\
&=\dashint_{B_{\rho}(y)}{}^t\!\Gamma_{kl}(\,\cdot\,,x)
=\dashint_{B_{\rho}(y)}\Gamma_{lk}(x,\,\cdot\,),
\end{split}
\end{equation}
i.e., we have the following representation for the averaged fundamental
matrix:
\begin{equation}
\label{eqE-68}
\vec{\Gamma}^\rho(x,y)=\dashint_{B_\rho(y)} \vec{\Gamma}(x,z)\,dz.
\end{equation}
Therefore, by the continuity, we obtain
\begin{equation}
\label{eqE-69}
\lim_{\rho\to 0} \vec{\Gamma}^\rho(x,y)=\vec{\Gamma}(x,y),
\quad x\neq y.
\end{equation}

\subsection{Properties of fundamental matrix}
We record what we obtained so far in the following theorem:
\begin{theorem} \label{thm:E-01}
Assume that operators $L$ and ${}^t\!L$ satisfy the property (H).
Then, there exists a unique fundamental matrix
$\vec{\Gamma}(x,y)=(\Gamma_{ij}(x,y))_{i,j=1}^N$ $(x\neq y)$
which is continuous in $\set{(x,y)\in\bR^n\times\bR^n:x\neq y}$
and such that $\vec\Gamma(x,\,\cdot\,)$ is locally integrable in
$\bR^n$ for all $x\in\bR^n$ and that for all
$\vec{f}=(f^1,\ldots,f^N)^T \in C^\infty_c(\bR^n)^N$,
the function $\vec{u}=(u^1,\ldots,u^N)^T$ given by
\begin{equation}
\label{eqE-70}
\vec{u}(x):=\int_{\bR^n} \vec{\Gamma}(x,y)\vec{f}(y)\,dy
\end{equation}
belongs to $Y^{1,2}(\bR^n)^N$ and satisfies $L\vec{u}=\vec{f}$ in the sense
\begin{equation}
\label{eqE-71}
\int_{\bR^n} A^{\alpha\beta}_{ij} D_\beta u^j D_\alpha \phi^i = \int_{\bR^n}
f^i \phi^i,\quad \forall \vec{\phi}\in C^\infty_c(\bR^n)^N.
\end{equation}
Moreover, $\vec{\Gamma}(x,y)$ has the property
\begin{equation}
\label{eqE-72}
\int_{\bR^n}A^{\alpha\beta}_{ij} D_\beta \Gamma_{jk}(\,\cdot\,,y)
D_\alpha \phi^i = \phi^k(y),
\quad\forall\vec{\phi}\in C^\infty_c(\bR^n)^N.
\end{equation}
Furthermore, $\vec{\Gamma}(x,y)$ satisfies the following estimates:
\begin{gather}
\label{eqE-73}
\norm{\vec{\Gamma}(\,\cdot\,,y)}_{Y^{1,2}(\bR^n\setminus B_r(y))}+
\norm{\vec{\Gamma}(x,\,\cdot\,)}_{Y^{1,2}(\bR^n\setminus B_r(x))}
\le C r^{1-\frac{n}{2}}, \quad \forall r>0,\\
\norm{\vec{\Gamma}(\,\cdot\,,y)}_{L^p(B_r(y))}+
\norm{\vec{\Gamma}(x,\,\cdot\,)}_{L^p(B_r(x))} \le C_p r^{2-n+\frac{n}{p}},
\quad \forall p\in [1,\tfrac{n}{n-2}),\\
\norm{D\vec{\Gamma}(\,\cdot\,,y)}_{L^p(B_r(y))}+
\norm{D\vec{\Gamma}(x,\,\cdot\,)}_{L^p(B_r(x))}\le C_p r^{1-n+\frac{n}{p}},
\quad\forall p\in [1,\tfrac{n}{n-1}),\\
\abs{\set{x\in\bR^n:\abs{\vec{\Gamma}(x,y)}>t}}+
\abs{\set{y\in\bR^n:\abs{\vec{\Gamma}(x,y)}>t}}\le
C t^{-\frac{n}{n-2}},\\
\abs{\set{x\in\bR^n:\abs{D_x \vec{\Gamma}(x,y)}>t}}+
\abs{\set{y\in\bR^n:\abs{D_y \vec{\Gamma}(x,y)}>t}}\le
C t^{-\frac{n}{n-1}},\\
\abs{\vec{\Gamma}(x,y)}\le C\abs{x-y}^{2-n},\quad \forall x\neq y,\\
\abs{\vec{\Gamma}(x,y)-\vec{\Gamma}(z,y)}\le 
C\abs{x-z}^{\mu_0} \abs{x-y}^{2-n-\mu_0}\quad\text{if }\abs{x-z}<\abs{x-y}/2,\\
\label{eqE-80}
\abs{\vec{\Gamma}(x,y)-\vec{\Gamma}(x,z)}\le 
C\abs{y-z}^{\mu_0} \abs{x-y}^{2-n-\mu_0}\quad\text{if }\abs{y-z}<\abs{x-y}/2,
\end{gather}
where $C=C(n,N,\lambda,\Lambda,\mu_0,H_0)>0$ and 
$C_p=C_p(n,N,\lambda,\Lambda,\mu_0,H_0,p)>0$.
\end{theorem}
\begin{proof}
Let $\vec{\Gamma}^\rho(x,y)$ and $\vec{\Gamma}(x,y)$ be constructed as above.
We have already seen that $\vec{\Gamma}$
is continuous in $\set{(x,y)\in \bR^n\times\bR^n:x\neq y}$ and satisfies
all the properties \eqref{eqE-72} -- \eqref{eqE-80}.
By using the Lax-Milgram lemma as in Sec.~\ref{sec:E-01},
we find that for all $\vec{f}\in C^\infty_c(\bR^n)^N$, there is a unique
$\vec{u}\in Y^{1,2}(\bR^n)^N$ satisfying
\begin{equation*}
\int_{\bR^n} A^{\alpha\beta}_{ij} D_\beta u^j D_\alpha v^i = \int_{\bR^n}
f^i v^i,\quad \forall \vec{v}\in Y^{1,2}(\bR^n)^N.
\end{equation*}
If we set $v^i=\Gamma^\rho_{ki}(x,\,\cdot\,)$ above, then
\eqref{eqE-07} together with \eqref{eqE-66} implies that
\begin{equation}
\label{eqE-81}
\int_{\bR^n} \Gamma^\rho_{ki}(x,\,\cdot\,) f^i=
\int_{\bR^n} {}^t\!A^{\beta\alpha}_{ji} D_\alpha {}^t\!\Gamma^\rho_{ik}
(\,\cdot\,,x) D_\beta u^j = \dashint_{B_\rho(x)} u^k.
\end{equation}
Assume that $\vec{f}$ is supported in $B_R(x)$ for some $R>0$.
Then, by \eqref{eqE-43} and \eqref{eqE-66} we have
\begin{equation*}
\lim_{\rho\to0}
\int_{\bR^n} \Gamma^\rho_{ki}(x,\,\cdot\,) f^i
= \lim_{\rho\to0}
\int_{B_R(x)} \Gamma^\rho_{ki}(x,\,\cdot\,) f^i
=\int_{\bR^n}\Gamma_{ki}(x,\,\cdot\,)f^i.
\end{equation*}
By the same argument which lead to \eqref{eqE-19}
in Section~\ref{sec:E-02}, we find that $\vec{u}$ is H\"older continuous.
Therefore, \eqref{eqE-70} follows by taking the limits in \eqref{eqE-81}.

Now, it only remains to prove the uniqueness.
Assume that $\tilde{\vec{\Gamma}}(x,y)$ is another matrix such that
$\tilde{\vec{\Gamma}}$ is continuous on
$\set{(x,y)\in \bR^n\times\bR^n:x\neq y}$
and such that $\tilde{\vec\Gamma}(x,\,\cdot\,)$ is locally integrable in $\bR^n$
for all $x\in\bR^n$ and that for all $\vec{f}\in C^\infty_c(\bR^n)^N$,
\begin{equation*}
\tilde{\vec{u}}(x):=\int_{\bR^n} \tilde{\vec{\Gamma}}(x,y)\vec{f}(y)\,dy
\end{equation*}
belongs to $Y^{1,2}(\bR^n)$ and satisfies $L\vec{u}=\vec{f}$ in the 
sense of \eqref{eqE-71}.
Then by the uniqueness in $\mathbf H=Y^{1,2}(\bR^n)^N$, we must have
$\vec{u}=\tilde{\vec{u}}$. Therefore, for all $x\in\bR^n$ we have
\begin{equation*}
\int_{\bR^n} (\vec{\Gamma}-\tilde{\vec{\Gamma}})(x,\cdot)\vec{f}=0,\quad
\forall \vec{f}\in C^\infty_c(\bR^n)^N,
\end{equation*}
and thus we have $\vec{\Gamma}\equiv\tilde{\vec{\Gamma}}$ in
$\set{(x,y)\in \bR^n\times\bR^n:x\neq y}$.
\end{proof}

\begin{theorem}
Assume that the operators $L$ and ${}^t\!L$ satisfy the property (H).
If $\vec{f}\in (L^{\frac{2n}{n+2}}(\bR^n)\cap L^p_{loc}(\bR^n))^N$
for some $p>n/2$, then there exists a unique
$\vec{u}$ in $Y^{1,2}(\bR^n)^N$ such that
\begin{equation}
\label{eqE-82}
\int_{\bR^n} A^{\alpha\beta}_{ij} D_\beta u^j D_\alpha v^i = \int_{\bR^n}
f^i v^i,\quad \forall \vec{v}\in Y^{1,2}(\bR^n)^N.
\end{equation}
Moreover, $\vec{u}$ is continuous and has the following representation:
\begin{equation}
\label{eqE-83}
u^k(x)=\int_{\bR^n} \Gamma_{ki}(x,y) f^i(y)\,dy,\quad
k=1,\ldots,N,
\end{equation}
where $\left(\Gamma_{ki}(x,y)\right)_{k,i=1}^N$ is
the fundamental matrix of $L$.
\end{theorem}
\begin{proof}
Since $\vec{f}\in L^{\frac{2n}{n+2}}(\bR^n)^N$, the same argument
as appears in Sec.~\ref{sec:E-01} implies that there is
$\vec{u}\in Y^{1,2}(\bR^n)^N$ satisfying \eqref{eqE-82}.
If we set $v^i=\Gamma^\rho_{ki}(x,\,\cdot\,)$ in \eqref{eqE-82}, then
\eqref{eqE-03} implies that
\begin{equation}
\label{eqE-84}
\int_{\bR^n} \Gamma^\rho_{ki}(x,\,\cdot\,) f^i=
\int_{\bR^n} {}^t\!A^{\beta\alpha}_{ji} D_\alpha {}^t\!\Gamma^\rho_{ik}
(\,\cdot\,,x) D_\beta u^j = \dashint_{B_\rho(x)} u^k.
\end{equation}
Next, note that \eqref{eqE-26}, \eqref{eqE-41}, and
the assumption $f^i\in L^{\frac{2n}{n+2}}(\bR^n)\cap L^p_{loc}(\bR^n)$
for some $p>n/2$, imply that
\begin{equation}
\label{eqE-85}
\begin{split}
\lim_{\rho\to0}
\int_{\bR^n} \Gamma^\rho_{ki}(x,\,\cdot\,) f^i
&= \lim_{\rho\to0}
\left(\int_{B_1(x)} \Gamma^\rho_{ki}(x,\,\cdot\,) f^i+
\int_{\bR^n\setminus B_1(x)}\Gamma^\rho_{ki}(x,\,\cdot\,)f^i\right)\\
&= \int_{B_1(x)} \Gamma_{ki}(x,\,\cdot\,) f^i+
\int_{\bR^n\setminus B_1(x)}\Gamma_{ki}(x,\,\cdot\,)f^i\\
&=\int_{\bR^n}\Gamma_{ki}(x,\,\cdot\,)f^i.
\end{split}
\end{equation}
Finally, by the same argument which lead to \eqref{eqE-19}
in Sec.~\ref{sec:E-02}, we find that $\vec{u}$ is H\"older continuous,
and thus \eqref{eqE-83} follows from \eqref{eqE-84} and \eqref{eqE-85}. 
\end{proof}

\begin{corollary}
Suppose that
$\vec{f}=(f^1,\ldots,f^N)^T$ has a bound 
\begin{equation}
\label{eqE-86}
\abs{\vec{f}(x)}\le C(1+\abs{x})^{-(1+n/2+\epsilon)}\quad \forall x\in\bR^n
\text{ for some $\epsilon>0$}.
\end{equation}
Then, $\vec{u}=(u^1,\ldots,u^N)^T$ given by \eqref{eqE-83}
is a unique $Y^{1,2}(\bR^n)^N$ solution of $L\vec{u}=\vec{f}$ in $\bR^n$ 
in the sense of \eqref{eqE-82}.
\end{corollary}
\begin{proof}
Note that \eqref{eqE-86} implies 
$\vec{f}\in (L^{\frac{2n}{n+2}}(\bR^n)\cap L^p(\bR^n))^N$.
\end{proof}

\begin{theorem}
Assume that $L$ and ${}^t\!L$ satisfy the property (H). If
$\vec{f}\in Y^{1,2}(\bR^n)^N$
satisfies $D\vec{f}\in L^p_{loc}(\bR^n)^{N\times n}$
for some $p>n$, then
\begin{equation}
\label{eqE-87}
f^k(x)=\int_{\bR^n} D_\alpha\Gamma_{ki}(x,\,\cdot\,)
A^{\alpha\beta}_{ij}D_\beta f^j,\quad k=1,\ldots,N,
\end{equation}
where $\left(\Gamma_{ki}(x,y)\right)_{k,i=1}^N$ is
the fundamental matrix of $L$.
\end{theorem}
\begin{proof}
We denote by ${}^t\!\vec{\Gamma}^\rho$ the averaged fundamental
matrix of ${}^t\!L$.
Recall that columns of ${}^t\!\vec{\Gamma}^\rho$ belong to $\mathbf H$.
Then, by \eqref{eqE-07} we have
\begin{equation*}
\int_{\bR^n} {}^t\!A^{\beta\alpha}_{ji}
D_\alpha{}^t\!\Gamma^\rho_{ik}(\,\cdot\,,x) D_\beta f^j
=\dashint_{B_\rho(x)}f^k.
\end{equation*}
As in \eqref{eqE-85},
the assumption $D\vec{f}\in L^p_{loc}(\bR^n)^N$ for $p>n$,
together with \eqref{eqE-34} and \eqref{eqE-37} yields
\begin{equation}
\label{eqE-89}
\begin{split}
\lim_{\rho\to0}\,\dashint_{B_\rho(x)}f^k
&= \lim_{\rho\to0}
\left(\int_{B_1(x)} + \int_{\bR^n\setminus B_1(x)}
{}^t\!A^{\beta\alpha}_{ji}
D_\alpha{}^t\!\Gamma^\rho_{ik}(\,\cdot\,,x) D_\beta f^j \right)\\
&= \int_{B_1(x)} + \int_{\bR^n\setminus B_1(x)}
{}^t\!A^{\beta\alpha}_{ji}
D_\alpha{}^t\!\Gamma_{ik}(\,\cdot\,,x) D_\beta f^j\\
&= \int_{\bR^n} A^{\alpha\beta}_{ij}
D_\alpha\Gamma_{ki}(x,\,\cdot\,) D_\beta f^j.
\end{split}
\end{equation}
where we used \eqref{eqE-66} in the last step.
By the Morrey's inequality \cite{Morrey},
$\vec{f}$ is continuous and thus \eqref{eqE-87} follows from \eqref{eqE-89}.
\end{proof}

\begin{corollary}
Assume that $L$, ${}^t\!L$, $\tilde L$, and ${}^t\!\tilde{L}$ satisfy
the property (H).
Denote by $\vec{\Gamma}$ and $\tilde{\vec{\Gamma}}$ the fundamental matrices
of $L$ and  $\tilde L$, respectively.
If the coefficients $A^{\alpha\beta}_{ij}$ of $L$
and $\tilde{A}^{\alpha\beta}_{ij}$ of $\tilde{L}$ are H\"older continuous,
then
\begin{equation}
\label{eqE-90}
\tilde{\Gamma}_{lm}(x,y)= \Gamma_{lm}(x,y)+
\int_{\bR^n} D_\alpha\Gamma_{li}(x,\,\cdot\,)
(A^{\alpha\beta}_{ij}- \tilde{A}^{\alpha\beta}_{ij})
D_\beta\tilde{\Gamma}{}_{jm}(\,\cdot\,,y),\quad x\neq y.
\end{equation}
\end{corollary}
\begin{proof}
We denote by $\vec{\Gamma}^\rho$ and $\tilde{\vec{\Gamma}}{}^\rho$
($\rho< \abs{x-y}/4$)
the averaged fundamental matrices of $L$ and $\tilde{L}$ respectively.
Recall that columns of $\vec{\Gamma}^\rho$ and $\tilde{\vec{\Gamma}}{}^\rho$
belong to $\mathbf H$.
Moreover, since we assume that the coefficients are H\"older continuous,
the standard elliptic theory, \eqref{eqE-54}, and \eqref{eqE-68}
implies that $D \vec{\Gamma}^\rho(x,\,\cdot\,)$ and
$D \tilde{\vec{\Gamma}}{}^\rho(\,\cdot\,,y)$ are
locally bounded.
Therefore, by setting $f^j= \tilde\Gamma{}^\rho_{jm}(\,\cdot\,,y)$
in \eqref{eqE-87} we have 
\begin{equation}
\label{eqE-91}
\tilde\Gamma{}^\rho_{lm}(x,y)= \int_{\bR^n}
D_\alpha\Gamma_{li}(x,\,\cdot\,) A^{\alpha\beta}_{ij}
D_\beta\tilde\Gamma{}^\rho_{jm}(\,\cdot\,,y),
\end{equation}
Next, set $f^j= \Gamma{}^\rho_{lj}(x,\,\cdot\,)$
and apply \eqref{eqE-87}
with $L$ replaced by ${}^t\!\tilde{L}$ to get
\begin{equation*}
\Gamma{}^\rho_{lm}(x,y)= \int_{\bR^n}
D_\alpha{}^t\!\tilde{\Gamma}_{mi}(y,\,\cdot\,)
{}^t\!\tilde{A}{}^{\alpha\beta}_{ij}
D_\beta \Gamma{}^\rho_{lj}(x,\,\cdot\,).
\end{equation*}
By using \eqref{eqE-66} and interchanging indices
($\alpha\leftrightarrow\beta$, $i\leftrightarrow j$), we obtain
\begin{equation}
\label{eqE-93}
\Gamma{}^\rho_{lm}(x,y)= \int_{\bR^n}
D_\alpha \Gamma{}^\rho_{li}(x,\,\cdot\,) \tilde{A}{}^{\alpha\beta}_{ij}
D_\beta\tilde{\Gamma}_{jm}(\,\cdot\,,y).
\end{equation}
Now, set $r=\abs{x-y}/4$ and split the integral \eqref{eqE-91}
into three pieces (recall $\rho <r$)
\begin{equation*}
\int_{B_r(x)}+ \int_{B_r(y)}+ \int_{\bR^n\setminus (B_r(x)\cup B_r(x))}
D_\alpha\Gamma_{li}(x,\,\cdot\,) A^{\alpha\beta}_{ij}
D_\beta\tilde\Gamma{}^\rho_{jm}(\,\cdot\,,y).
\end{equation*}
Since we assume that the coefficients are H\"older continuous,
it follows from the standard elliptic theory that
$D \vec{\Gamma}(x,\,\cdot\,)$ and $D \tilde{\vec{\Gamma}}(\,\cdot\,,y)$
are continuous (and thus bounded) on $B_r(y)$ and $B_r(x)$ respectively.
Moreover, \eqref{eqE-68} implies
\begin{equation*}
D\tilde{\vec{\Gamma}}{}^\rho(\,\cdot\,,y)\to
D\tilde{\vec{\Gamma}}{}(\,\cdot\,,y)\quad\text{uniformly on }
B_r(x)\text{ as }\rho\to 0.
\end{equation*}
Therefore, as in \eqref{eqE-89}, we may take the limit $\rho\to0$
in \eqref{eqE-91} to get
\begin{equation*}
\tilde\Gamma_{lm}(x,y)= \int_{\bR^n}
D_\alpha\Gamma_{li}(x,\,\cdot\,) A^{\alpha\beta}_{ij}
D_\beta\tilde\Gamma_{jm}(\,\cdot\,,y),
\end{equation*}
Similarly, by taking the limit $\rho\to0$ in \eqref{eqE-93}, we obtain
\begin{equation*}
\Gamma_{lm}(x,y)= \int_{\bR^n}
D_\alpha \Gamma_{li}(x,\,\cdot\,) \tilde{A}{}^{\alpha\beta}_{ij}
D_\beta\tilde{\Gamma}_{jm}(\,\cdot\,,y).
\end{equation*}
The proof is complete.
\end{proof}

\begin{remark}
We note that in terms of matrix multiplication
\eqref{eqE-83} is written as 
\begin{equation*}
\vec{u}(x)=\int_{\bR^n}\vec{\Gamma}(x,y)\vec{f}(y)\,dy,
\end{equation*}
where both $\vec{u}, \vec{f}$ are understood as column vectors.
Also, \eqref{eqE-90} reads
\begin{equation*}
\tilde{\vec{\Gamma}}(x,y)= \vec{\Gamma}(x,y)+
\int_{\bR^n} D_\alpha\vec{\Gamma}(x,\,\cdot\,)
(\vec{A}^{\alpha\beta}-\tilde{\vec{A}}{}^{\alpha\beta})
D_\beta\vec{\Gamma}(\,\cdot\,,y).
\end{equation*}
\end{remark}

\section{Green's matrix in general domains}\label{sec:G}
\subsection{Construction of Green's matrix}\label{sec:G-01}
In this section, we shall construct the Green's matrix in any open,  
connected set $\Omega\subset \bR^n$, where $n\ge 3$.
To construct the Green's matrix in $\Omega$, we need to adjust arguments
in Section~\ref{sec:E}.

Henceforth, we shall denote $\Omega_r(y):=\Omega\cap B_r(y)$
and $d_y:=\dist(y,\partial\Omega)$.
Also, as in Section~\ref{sec:E}, we use the letter $C$ to denote
a constant depending on $n$, $N$, $\lambda$, $\Lambda$, $\mu_0$, $H_0$,
and sometimes on an exponent $p$ characterizing Lebesgue classes.

It is routine to check that for any given $y\in\Omega$ and $1\le k\le N$,
the linear functional $\vec{u}\mapsto \dashint_{\Omega_\rho(y)} u^k$
is bounded on $\mathbf H=Y^{1,2}_0(\Omega)^N$.
Therefore, by Lax-Milgram lemma, there exists a unique
$\vec{v}_\rho=\vec{v}_{\rho;y,k}\in \mathbf H$ such that
\begin{equation}
\label{eqG-01}
\int_{\Omega}A^{\alpha\beta}_{ij}D_\beta v^j_{\rho}\,D_\alpha u^i=
\dashint_{\Omega_\rho(y)}u^k,\quad\forall\vec{u}\in \mathbf H.
\end{equation}
Note that as in \eqref{eqE-05}, we have
\begin{equation}
\label{eqG-02}
\norm{D\vec{v}_\rho}_{L^2(\Omega)}=\norm{\vec{v}_\rho}_{\mathbf H}\le
C\,\Abs{\Omega_\rho(y)}^{\frac{2-n}{2n}}.
\end{equation}
We define the ``averaged Green's matrix''
$\vec{G}^\rho(\,\cdot\,,y)=(G_{jk}^\rho(\,\cdot\,,y))_{j,k=1}^N$ by
\begin{equation*}
G^\rho_{jk}(\,\cdot\,,y)=v^j_{\rho}=v^j_{\rho;y,k}.
\end{equation*}
Note that as in \eqref{eqE-07}, we have
\begin{equation}
\label{eqG-04}
\int_{\Omega}A^{\alpha\beta}_{ij} D_\beta G^\rho_{jk}(\,\cdot\,,y)
D_\alpha u^i = \dashint_{\Omega_\rho(y)}u^k,
\quad\forall\vec{u}\in {\mathbf H}.
\end{equation}

Next, observe that as in \eqref{eqE-09}--\eqref{eqE-12},
for any given $\vec{f}\in L^\infty_c(\Omega)^N$,
there exists a unique $\vec{u}\in {\mathbf H}$ such that
\begin{equation*}
\int_{\Omega} G^\rho_{ik}(\,\cdot\,,y) f^i=
\dashint_{\Omega_\rho(y)}u^k.
\end{equation*}
Moreover, as in \eqref{eqE-13}, we have
\begin{equation*}
\norm{D\vec{u}}_{L^2(\Omega)}\le C
\norm{\vec{f}}_{L^{2n/(n+2)}(\Omega)}.
\end{equation*}
Also, by following the argument as appears in Section~\ref{sec:E-02},
we find that if $\vec{f}$ is supported in $B_R(y)$, then we have
\begin{equation*}
\norm{\vec{u}}_{L^\infty(B_{R/4}(y))}
\le C R^{2-n/p} \norm{\vec{f}}_{L^p(B_R(y))},
\quad \forall R < d_y,
\quad \forall p>n/2.
\end{equation*}
Therefore, as in \eqref{eqE-22}, for any $\vec{f}\in L^\infty_c(B_R(y))$,
$R< d_y$, we have
\begin{equation*}
\abs{\int_{B_R(y)} G^\rho_{ik}(\,\cdot\,,y) f^i\,}\le
C R^{2-n/p} \norm{\vec{f}}_{L^p(B_R(y))},\quad
\forall \rho<R/4,\quad \forall p>n/2.
\end{equation*}
Therefore, as in \eqref{eqE-23},
we see that if $R< d_y$, then
\begin{equation*}
\norm{\vec{G}^\rho(\,\cdot\,,y)}_{L^q(B_R(y))}\le C R^{2-n+n/q},
\quad \forall \rho<R/4, \quad\forall q\in [1,\tfrac{n}{n-2}).
\end{equation*}
Then, by following the lines in \eqref{eqE-24}--\eqref{eqE-25},
we obtain
\begin{equation*}
\abs{\vec{G}^\rho(x,y)}\le C \abs{x-y}^{2-n}
\quad \text{if } \abs{x-y}<d_y/2,\quad
\forall \rho<\abs{x-y}/3.
\end{equation*}

Next, we shall derive an estimate corresponding to \eqref{eqE-34}.
Let $\eta\in C^\infty(\bR^n)$ be  a cut-off function such that
$0\le \eta\le 1$, $\eta\equiv 1$ outside $B_{R/2}(y)$,
$\eta\equiv 0$ on $B_{R/4}(y)$, and $\abs{D\eta}\le C/R$, where $R\le d_y$.
By setting $\vec{u}=\eta^2\vec{v}_\rho\in {\mathbf H}$ in \eqref{eqG-01}, we obtain
\begin{equation}
\label{eqG-11}
\begin{split}
\int_\Omega \eta^2\abs{D\vec{v}_\rho}^2
&\le C \int_\Omega \abs{D\eta}^2\abs{\vec{v}_\rho}^2
\le C R^{-2} \int_{B_{R/2}(y)\setminus B_{R/4}(y)}\abs{\vec{v}_\rho}^2\\
&\le C R^{-2} \int_{B_{R/2}(y)\setminus B_{R/4}(y)}\abs{x-y}^{2(2-n)}\,dx\\
&=C R^{-2} R^{4-n}= C R^{2-n},\quad \forall \rho< R/12.
\end{split}
\end{equation}
Therefore, we have ($r=R/2$)
\begin{equation}
\label{eqG-12}
\int_{\Omega\setminus B_r(y)} \abs{D\vec{G}^\rho(\,\cdot\,,y)}^2 \le C r^{2-n},
\quad \forall \rho< r/6,
\quad \forall r< d_y/2.
\end{equation}
On the other hand, \eqref{eqG-02} implies that
if $\rho\ge r/6$, then
\begin{equation}
\label{eqG-13}
\int_{\Omega\setminus B_r(y)}\abs{D\vec{G}^\rho(\,\cdot\,,y)}^2
\le \int_{\Omega}\abs{D\vec{G}^\rho(\,\cdot\,,y)}^2
\le C\,\Abs{\Omega_\rho(y)}^{\frac{2-n}{n}}
\le C r^{2-n}.
\end{equation}
Therefore, by combining \eqref{eqG-12} and \eqref{eqG-13}, we obtain
\begin{equation}
\label{eqG-14}
\int_{\Omega\setminus B_r(y)} \abs{D\vec{G}^\rho(\,\cdot\,,y)}^2 \le C r^{2-n},
\quad \forall r< d_y/2,\quad \forall \rho>0.
\end{equation}

From the estimate \eqref{eqG-14}, which corresponds
to \eqref{eqE-34}, we can derive an estimate corresponding
to \eqref{eqE-37} as follows.
By following the lines between \eqref{eqE-34} and \eqref{eqE-36},
we obtain
\begin{equation}
\label{eqG-15}
\abs{\set{x\in\Omega:\abs{D_x \vec{G}^\rho(x,y)}>t}}\le
C t^{-\frac{n}{n-1}},\quad \forall \rho>0\quad\text{if }
t> (d_y/2)^{1-n}.
\end{equation}
Then, by following lines \eqref{eqE-37}--\eqref{eqE-40},
we find (set $\tau= (R/2)^{1-n}$)
\begin{equation}
\label{eqG-16}
\int_{B_R(y)} \abs{D\vec{G}^\rho(\,\cdot\,,y)}^p \le C R^{p(1-n)+n},
\quad \forall R< d_y,\quad\forall \rho>0,
\quad \forall p\in (0,\tfrac{n}{n-1}).
\end{equation}

Now, we will derive estimates corresponding \eqref{eqE-26}
and \eqref{eqE-41}.
Let $\eta$ be the same as in \eqref{eqG-11}.
Note that \eqref{eqG-11} and \eqref{eqG-14} implies that for $R< d_y$,
\begin{equation}
\label{eqG-17}
\int_\Omega \abs{D(\eta\vec{v}_\rho)}^2
\le 2\int_\Omega \eta^2 \abs{D\vec{v}_\rho}^2
+2\int_\Omega \abs{D\eta}^2 \abs{\vec{v}_\rho}^2
\le C R^{2-n},\quad\forall \rho<R/12.
\end{equation}
Since $\eta\vec{v}_\rho\in {\mathbf H}=Y^{1,2}_0(\Omega)$,
it follows from \eqref{eqG-17} and \eqref{eqP-14} that
\begin{equation}
\label{eqG-18}
\int_{\Omega\setminus B_r(y)} \abs{\vec{v}_\rho}{}^{2^*}
\le C r^{-n}, \quad\forall r< d_y/2, \quad\forall \rho<r/6.
\end{equation}
On the other hand, if $\rho\ge r/6$, then \eqref{eqG-02} implies
\begin{equation}
\label{eqG-19}
\begin{split}
\int_{\Omega\setminus B_r(y)} \abs{\vec{v}_\rho}{}^{2^*}
&\le
\int_{\Omega} \abs{\vec{v}_\rho}{}^{2^*}
\le C \left(\int_{\Omega} \abs{D\vec{v}_\rho}{}^2\right)^{2^*/2}\\
&\le C \abs{\Omega_\rho}^{-1}\le C r^{-n}.
\end{split}
\end{equation}
Therefore, by combining \eqref{eqG-18} and \eqref{eqG-19}, we obtain
\begin{equation}
\label{eqG-20}
\int_{\Omega\setminus B_r(y)} \abs{\vec{G}^\rho(\,\cdot\,,y)}{}^{2^*}
\le C r^{-n},
\quad \forall r<d_y/2,\quad \forall \rho>0.
\end{equation}
As in Section~\ref{sec:E-03},
the above estimate \eqref{eqG-20} yields
\begin{equation}
\label{eqG-21}
\abs{\set{x\in\Omega:\abs{\vec{G}^\rho(x,y)}>t}}\le
C t^{-\frac{n}{n-2}},\quad \forall \rho>0\quad\text{if }
t> (d_y/2)^{2-n}.
\end{equation}
Then, as we argued in \eqref{eqG-16}, we find (set $\tau= (R/2)^{2-n}$)
\begin{equation}
\label{eqG-22}
\int_{B_R(y)} \abs{\vec{G}^\rho(\,\cdot\,,y)}^p \le C R^{p(2-n)+n},
\quad \forall R< d_y,\quad\forall \rho>0,
\quad \forall p\in (0,\tfrac{n}{n-2}).
\end{equation}

Now, observe that \eqref{eqG-16} and \eqref{eqG-22} in particular
imply that
\begin{equation}
\label{eqG-23}
\norm{\vec{G}^\rho(\,\cdot\,,y)}_{W^{1,p}(B_{d_y}(y))}\le
C(d_y)\text{  for some }p\in (1,\tfrac{n}{n-1}),\text{ uniformly in }\rho.
\end{equation}
Therefore, from \eqref{eqG-23} together with \eqref{eqG-14} and
\eqref{eqG-20}, it follows that there exist a sequence
$\set{\rho_\mu}_{\mu=1}^\infty$ tending to $0$ and functions
$\vec{G}(\,\cdot\,,y)$ and $\tilde{\vec{G}}(\,\cdot\,,y)$ such that
\begin{gather}
\label{eqG-24}
\vec{G}^{\rho_\mu}(\,\cdot\,,y)\rightharpoonup \vec{G}(\,\cdot\,,y)
\quad\text{in}\quad W^{1,p}(B_{d_y}(y))^{N\times N}\quad\text{and}\\
\label{eqG-25}
\vec{G}^{\rho_\mu}(\,\cdot\,,y)\rightharpoonup \tilde{\vec{G}}(\,\cdot\,,y)
\quad\text{in}\quad Y^{1,2}(\Omega\setminus B_{d_y/2}(y))^{N\times N}
\text{ as } \mu\to\infty.
\end{gather}
Since $\vec{G}(\,\cdot\,,y)\equiv\tilde{\vec{G}}(\,\cdot\,,y)$
on $B_{d_y}(y)\setminus B_{d_y/2}(y)$,
we shall extend $\vec{G}(\,\cdot\,,y)$ to entire $\Omega$
by setting $\vec{G}(\,\cdot\,,y)=\tilde{\vec{G}}(\,\cdot\,,y)$
on $\Omega\setminus B_{d_y}(y)$ but still call it 
$\vec{G}(\,\cdot\,,y)$ in the sequel.
Moreover, by applying a diagonalization process
and passing to a subsequence, if necessary, we may assume that
\begin{equation}
\label{eqG-26}
\vec{G}^{\rho_{\mu}}(\,\cdot\,,y)\rightharpoonup \vec{G}(\,\cdot\,,y)
\quad\text{in }Y^{1,2}(\Omega\setminus B_r(y))^{N\times N}\text{ as }
\mu\to\infty,\quad\forall r< d_y.
\end{equation}

We claim that the following holds:
\begin{equation}
\label{eqG-27}
\int_{\Omega}A^{\alpha\beta}_{ij} D_\beta G_{jk}(\,\cdot\,,y)
D_\alpha \phi^i = \phi^k(y),\quad
\forall \vec{\phi}\in C^\infty_c(\Omega)^N.
\end{equation}
To see \eqref{eqG-27}, write $\vec{\phi}= \eta \vec{\phi}+ (1-\eta)
\vec{\phi}$, where $\eta\in C^\infty_c(B_{d_y}(y))$ is a cut-off
function satisfying $\eta\equiv 1$ on $B_{d_y/2}(y)$.
Then, \eqref{eqG-04}, \eqref{eqG-24}, and \eqref{eqG-26} yield
\begin{equation*}
\begin{split}
\phi^k(y)
&=\lim_{\mu\to\infty} \dashint_{\Omega_{\rho_\mu}(y)} \eta \phi^k +
\lim_{\mu\to\infty} \dashint_{\Omega_{\rho_\mu}(y)}(1-\eta) \phi^k\\
&=\lim_{\mu\to\infty} \int_{\Omega}A^{\alpha\beta}_{ij}
D_\beta G^{\rho_\mu}_{jk}(\cdot,y) D_\alpha (\eta \phi^i)
+ \lim_{\mu\to\infty}\int_{\Omega}A^{\alpha\beta}_{ij}
D_\beta G^{\rho_\mu}_{jk}(\cdot,y) D_\alpha ((1-\eta) \phi^i)\\
&= \int_{\Omega}A^{\alpha\beta}_{ij} D_\beta G_{jk}(\,\cdot\,,y)
D_\alpha (\eta \phi^i)
+ \int_{\Omega}A^{\alpha\beta}_{ij} D_\beta G_{jk}(\,\cdot\,,y)
D_\alpha ((1-\eta) \phi^i)\\
&=\int_{\Omega}A^{\alpha\beta}_{ij}D_\beta G_{jk}(\,\cdot\,,y)
D_\alpha \phi^i\quad\text{as desired}.
\end{split}
\end{equation*}
Next, we claim that $\vec{G}(\,\cdot\,,y)=0$ on $\partial\Omega$
in the sense that
for all $\eta\in C^\infty_c(\Omega)$ satisfying
$\eta\equiv 1$ on $B_r(y)$ for some $r<d_y$,
we have
\begin{equation*}
(1-\eta)\vec{G}(\,\cdot\,,y)\in Y^{1,2}_0(\Omega)^{N\times N}.
\end{equation*}
To see this, it is enough to show that
\begin{equation}
\label{eqG-29}
(1-\eta)\vec{G}^{\rho_\mu}(\,\cdot\,,y)\rightharpoonup
(1-\eta)\vec{G}(\,\cdot\,,y)
\quad\text{in }Y^{1,2}(\Omega)^{N\times N}\text{ as }\mu\to\infty,
\end{equation}
for $(1-\eta)\vec{G}^{\rho_\mu}(\,\cdot\,,y)\in Y^{1,2}_0(\Omega)^{N\times N}$
for all $\mu\ge 1$ and $Y^{1,2}_0(\Omega)$ is weakly closed in
$Y^{1,2}(\Omega)$ by Mazur's theorem.
To show \eqref{eqG-29}, we note that \eqref{eqG-26} yields
\begin{equation*}
\begin{split}
\int_\Omega (1-\eta) G_{kl}(\,\cdot\,,y)\phi&=
\int_\Omega G_{kl}(\,\cdot\,,y)(1-\eta) \phi=
\lim_{\mu\to\infty}\int_\Omega G_{kl}^{\rho_\mu}(\,\cdot\,,y)(1-\eta)\phi\\
&=\lim_{\mu\to\infty}\int_\Omega (1-\eta) G_{kl}^{\rho_\mu}(\,\cdot\,,y)\phi,
\quad\forall\phi\in L^{\frac{2n}{n+2}}(\Omega),
\end{split}
\end{equation*}
\begin{equation*}
\begin{split}
\int_\Omega D((1-\eta)&G_{kl}(\,\cdot\,,y))\cdot\vec{\psi}
=-\int_\Omega G_{kl}(\,\cdot\,,y) D\eta\cdot\vec{\psi}+
\int_\Omega DG_{kl}(\,\cdot\,,y)\cdot (1-\eta)\vec{\psi}\\
&=-\lim_{\mu\to\infty}
\int_\Omega G_{kl}^{\rho_\mu}(\,\cdot\,,y) D\eta\cdot\vec{\psi}
+\lim_{\mu\to\infty}
\int_\Omega DG_{kl}^{\rho_\mu}(\,\cdot\,,y)\cdot (1-\eta)\vec{\psi}\\
&=\lim_{\mu\to\infty}\int_\Omega
D((1-\eta) G_{kl}^{\rho_\mu}(\,\cdot\,,y))\cdot \vec{\psi},
\quad \forall \vec{\psi}\in L^2(\Omega)^N.
\end{split}
\end{equation*}

By using the same duality argument as in \eqref{eqE-45},
we derive the following estimates that correspond to
\eqref{eqE-46}--\eqref{eqE-52}:
\begin{gather}
\label{eqG-32}
\norm{\vec{G}(\,\cdot\,,y)}_{L^p(B_r(y))}\le C_p\,r^{2-n+n/p},
\quad \forall r< d_y,
\quad \forall p\in [1,\tfrac{n}{n-2}),\\
\norm{D\vec{G}(\,\cdot\,,y)}_{L^p(B_r(y))}\le C_p\,r^{1-n+n/p},
\quad \forall r< d_y,
\quad \forall p\in [1,\tfrac{n}{n-1}),\\
\label{eqG-34}
\norm{\vec{G}(\,\cdot\,,y)}_{Y^{1,2}(\Omega\setminus B_r(y))}
\le C r^{1-n/2},\quad \forall r< d_y/2,\\
\label{eqG-35}
\abs{\set{x\in\Omega:\abs{\vec{G}(x,y)}>t}}\le C t^{-\frac{n}{n-2}},
\quad\forall t> (d_y/2)^{2-n},\\
\label{eqG-36}
\abs{\set{x\in\Omega:\abs{D_x \vec{G}(x,y)}>t}}\le C t^{-\frac{n}{n-1}},
\quad\forall t> (d_y/2)^{1-n}.
\end{gather}

Also, we obtain pointwise bound and H\"older continuity estimate for
$\vec{G}(\,\cdot\,,y)$ corresponding to \eqref{eqE-54} and \eqref{eqE-55},
respectively, as follows.
Denote by $\vec{v}^T$ the $k$-th column of $\vec{G}(\,\cdot\,,y)$
and set $R:=\bar{d}_{x,y}/2$, where
\begin{equation}
\label{eqG-37}
\bar{d}_{x,y}:=\min(d_x, d_y,\abs{x-y}).
\end{equation}
Since $\vec{v}$ is a weak solution of $L\vec{u}=0$ in
$B_{3R/2}(x)\subset\Omega\setminus B_{R/2}(y)$, it follows from
\eqref{eqP-11} and \eqref{eqG-34} that
\begin{equation*}
\abs{\vec{v}(x)}
\le C R^{(2-n)/2}\norm{\vec{v}}_{L^{2^*}(\Omega\setminus B_{R/2}(y))}
\le C R^{2-n},
\end{equation*}
which in turn implies that
\begin{equation}
\label{eqG-39}
\abs{\vec{G}(x,y)}\le C \bar{d}_{x,y}^{2-n},\quad
\text{where }\bar{d}_{x,y}:=\min(d_x, d_y,\abs{x-y}).
\end{equation}
In particular, we have
\begin{equation}
\label{eqG-40}
\abs{\vec{G}(x,y)}\le C \abs{x-y}^{2-n}
\quad\text{if }\abs{x-y}<d_x/2
\text{ or } \abs{x-y}<d_y/2.
\end{equation}
Similarly, it follows from \eqref{eqP-12} and \eqref{eqG-34} that
\begin{equation}
\label{eqG-41}
[\vec{v}]^2_{C^{\mu_0}(B_R(x))}
\le C R^{2-n-2\mu_0} \int_{B_{3R/2}(x)} \abs{D\vec{v}}^2
\le C R^{2(2-n-\mu_0)}.
\end{equation}
Therefore, we find that
\begin{equation}
\label{eqG-42}
\abs{\vec{G}(x,y)-\vec{G}(z,y)}\le
C \abs{x-z}^{\mu_0} \bar{d}_{x,y}^{2-n-\mu_0}\quad\text{if }
\abs{x-z}<\bar{d}_{x,y}/2,
\end{equation}
where $\bar{d}_{x,y}$ is given by \eqref{eqG-37}.

Denote by ${}^t\!\vec{G}^\sigma(\,\cdot\,,x)$ the averaged
Green's matrix of ${}^t\!L$ in $\Omega$ with a pole at $x\in\Omega$.
Observe that we have an identity corresponding to \eqref{eqE-59}.
\begin{equation}
\label{eqG-43}
\dashint_{\Omega_{\rho}(y)}{}^t\!G_{kl}^{\sigma}(\,\cdot\,,x)
=\dashint_{\Omega_{\sigma}(x)}G_{lk}^{\rho}(\,\cdot\,,y).
\end{equation}
Let ${}^t\!\vec{G}(\,\cdot\,,x)$ be a Green's matrix of ${}^t\!L$ in
$\Omega$ with a pole at $x\in\Omega$ that is obtained by a sequence
$\set{\sigma_\nu}_{\nu=1}^\infty$ tending to $0$.
Then, by a similar argument as appears in Section~\ref{sec:E-06}, we obtain
\begin{equation}
G_{lk}(x,y)={}^t\!G_{kl}(y,x),
\quad \forall k,l=1,\ldots,N,\quad \forall x,y \in\Omega,\quad x\neq y,
\end{equation}
which is equivalent to say
\begin{equation}
\label{eqG-45}
\vec{G}(x,y)={}^t\!\vec{G}(y,x)^T,\quad \forall x,y\in\Omega,\quad x\neq y.
\end{equation}
Using \eqref{eqG-45}, we find that $\vec{G}(x,\,\cdot\,)$ satisfies
the estimates corresponding to
\eqref{eqG-32}--\eqref{eqG-36} and \eqref{eqG-42}.
Moreover, by following the lines \eqref{eqE-67}--\eqref{eqE-68}
and using \eqref{eqG-43} we obtain 
\begin{equation}
\label{eqG-52}
\vec{G}^\rho(x,y)=\dashint_{\Omega_\rho(y)} \vec{G}(x,z)\,dz.
\end{equation}
Therefore, by the continuity, we find
\begin{equation}
\lim_{\rho\to 0} \vec{G}^\rho(x,y)=\vec{G}(x,y), \quad \forall x,y\in\Omega,
\quad x\neq y.
\end{equation}

Finally, we summarize what we obtained so far in the following theorem.
\begin{theorem}   \label{thm:G-01}
Let $\Omega$ be an open connected set in $\bR^n$. Denote
$d_x:=\dist(x,\partial\Omega)$ for $x\in\Omega$; we set $d_x=\infty$ if 
$\Omega=\bR^n$.
Assume that operators $L$ and ${}^t\!L$ satisfy the property (H).
Then, there exists a unique Green's matrix
$\vec{G}(x,y)=(G_{ij}(x,y))_{i,j=1}^N$ ($x,y\in\Omega, x\neq y$)
which is continuous in $\set{(x,y)\in\Omega\times\Omega:x\neq y}$
and such that $\vec G(x,\,\cdot\,)$ is locally integrable in
$\Omega$ for all $x\in\Omega$ and that for all
$\vec{f}=(f^1,\ldots,f^N)^T \in C^\infty_c(\Omega)^N$,
the function $\vec{u}=(u^1,\ldots,u^N)^T$ given by
\begin{equation}
\label{eqZ-70}
\vec{u}(x):=\int_{\Omega} \vec{G}(x,y)\vec{f}(y)\,dy
\end{equation}
belongs to $Y^{1,2}_0(\Omega)^N$ and satisfies $L\vec{u}=\vec{f}$ in the sense
\begin{equation}
\label{eqZ-71}
\int_{\Omega} A^{\alpha\beta}_{ij} D_\beta u^j D_\alpha \phi^i = \int_{\Omega}
f^i \phi^i,\quad \forall \vec{\phi}\in C^\infty_c(\Omega)^N.
\end{equation}
Moreover, $\vec{G}(x,y)$ has the properties that
\begin{equation}
\label{eqG-54}
\int_{\Omega}A^{\alpha\beta}_{ij} D_\beta G_{jk}(\,\cdot\,,y)
D_\alpha \phi^i = \phi^k(y),
\quad \forall \vec{\phi}\in C^\infty_c(\Omega)^N
\end{equation}
and that for all $\eta\in C^\infty_c(\Omega)$ satisfying
$\eta\equiv 1$ on $B_r(y)$ for some $r<d_y$,
\begin{equation}
(1-\eta)\vec{G}(\,\cdot\,,y)\in Y^{1,2}_0(\Omega)^{N\times N}.
\end{equation}
Furthermore, $\vec{G}(x,y)$ satisfies the following estimates:
\begin{gather}
\norm{\vec{G}(\,\cdot\,,y)}_{L^p(B_r(y))}\le C_p\,r^{2-n+n/p},
\quad \forall r< d_y, \quad \forall p\in [1,\tfrac{n}{n-2}),\\
\norm{\vec{G}(x,\,\cdot\,)}_{L^p(B_r(x))}\le C_p\,r^{2-n+n/p},
\quad \forall r< d_x, \quad \forall p\in [1,\tfrac{n}{n-2}),
\end{gather}
\begin{gather}
\norm{D\vec{G}(\,\cdot\,,y)}_{L^p(B_r(y))}\le C_p\,r^{1-n+n/p},
\quad \forall r< d_y, \quad \forall p\in [1,\tfrac{n}{n-1}),\\
\norm{D\vec{G}(x,\,\cdot\,)}_{L^p(B_r(x))}\le C_p\,r^{1-n+n/p},
\quad \forall r< d_x, \quad \forall p\in [1,\tfrac{n}{n-1}),
\end{gather}
\begin{gather}
\norm{\vec{G}(\,\cdot\,,y)}_{Y^{1,2}(\Omega\setminus B_r(y))}
\le C r^{1-n/2},\quad \forall r< d_y/2,\\
\norm{\vec{G}(x,\,\cdot\,)}_{Y^{1,2}(\Omega\setminus B_r(x))}
\le C r^{1-n/2},\quad \forall r< d_x/2,
\end{gather}
\begin{gather}
\abs{\set{x\in\Omega:\abs{\vec{G}(x,y)}>t}}\le C t^{-\frac{n}{n-2}},
\quad\forall t> (d_y/2)^{2-n},\\
\abs{\set{y\in\Omega:\abs{\vec{G}(x,y)}>t}}\le C t^{-\frac{n}{n-2}},
\quad\forall t> (d_x/2)^{2-n},
\end{gather}
\begin{gather}
\abs{\set{x\in\Omega:\abs{D_x \vec{G}(x,y)}>t}}\le C t^{-\frac{n}{n-1}},
\quad\forall t> (d_y/2)^{1-n}.\\
\abs{\set{y\in\Omega:\abs{D_y \vec{G}(x,y)}>t}}\le C t^{-\frac{n}{n-1}},
\quad\forall t> (d_x/2)^{1-n},
\end{gather}
\begin{equation}
\abs{\vec{G}(x,y)}\le C \bar{d}_{x,y}^{2-n},\quad
\text{where }\bar{d}_{x,y}:=\min(d_x, d_y,\abs{x-y}),
\end{equation}
\begin{gather}
\abs{\vec{G}(x,y)-\vec{G}(z,y)}\le
C \abs{x-z}^{\mu_0} \bar{d}_{x,y}^{2-n-\mu_0}\quad\text{if }
\abs{x-z}<\bar{d}_{x,y}/2,\\
\label{eqG-55}
\abs{\vec{G}(x,y)-\vec{G}(x,z)}\le
C \abs{y-z}^{\mu_0} \bar{d}_{x,y}^{2-n-\mu_0}\quad\text{if }
\abs{y-z}<\bar{d}_{x,y}/2,
\end{gather}
where $C=C(n,N,\lambda,\Lambda,\mu_0,H_0)>0$ and 
$C_p=C_p(n,N,\lambda,\Lambda,\mu_0,H_0,p)>0$.
\end{theorem}
\begin{proof}
Let $\vec{G}^\rho(x,y)$ and $\vec{G}(x,y)$ be constructed as above.
We have already seen that $\vec{G}$
is continuous on $\set{(x,y)\in \Omega\times\Omega:x\neq y}$ and satisfies
all the properties \eqref{eqG-54} -- \eqref{eqG-55}.
Also, as in the proof of Theorem~\ref{thm:E-01},
we find that for all $\vec{f}\in C^\infty_c(\Omega)^N$, there is a unique
$\vec{u}\in (Y^{1,2}_0(\Omega)\cap C(\Omega))^N$ satisfying
\begin{equation*}
\int_{\Omega} A^{\alpha\beta}_{ij} D_\beta u^j D_\alpha v^i = \int_{\Omega}
f^i v^i,\quad \forall \vec{v}\in Y^{1,2}_0(\Omega)^N.
\end{equation*}
If we set $v^i=G^\rho_{ki}(x,\,\cdot\,)$ above, then
by \eqref{eqG-04} and \eqref{eqG-45}, we find
\begin{equation}
\label{eqZ-81}
\int_{\Omega} G^\rho_{ki}(x,\,\cdot\,) f^i=
\int_{\Omega} {}^t\!A^{\beta\alpha}_{ji} D_\alpha {}^t\!G^\rho_{ik}
(\,\cdot\,,x) D_\beta u^j = \dashint_{\Omega_\rho(x)} u^k.
\end{equation}
Fix $r<d_x/2$.
By \eqref{eqG-24}, \eqref{eqG-25}, and \eqref{eqG-45},
we have
\begin{equation*}
\begin{split}
\lim_{\mu\to\infty}
\int_{\Omega} G^{\rho_\mu}_{ki}(x,\,\cdot\,) f^i
&= \lim_{\mu\to\infty}
\left(\int_{B_r(x)} G^{\rho_\mu}_{ki}(x,\,\cdot\,) f^i+
\int_{\Omega\setminus B_r(x)}G^{\rho_\mu}_{ki}(x,\,\cdot\,)f^i\right)\\
&= \int_{B_1(x)} G_{ki}(x,\,\cdot\,) f^i+
\int_{\Omega\setminus B_1(x)}G_{ki}(x,\,\cdot\,)f^i\\
&=\int_{\Omega}G_{ki}(x,\,\cdot\,)f^i.
\end{split}
\end{equation*}
Therefore, \eqref{eqZ-70} follows by taking the limits in \eqref{eqZ-81}.
By proceeding as in the proof of Theorem~\ref{thm:E-01}, we also derive
the uniqueness of Green's matrix in $\Omega$.
\end{proof}

\subsection{Boundary regularity}\label{sec:G-02}
Let $\Sigma$ be any subset of $\overline{\Omega}$  and 
$u$ be a $W^{1,2}(\Omega)$ function. Then we shall say $u=0$ on
$\Sigma$ (in the sense of $W^{1,2}(\Omega)$) 
if $u$ is a limit in $W^{1,2}(\Omega)$ of a sequence of functions
in $C^\infty_c(\overline{\Omega}\setminus\Sigma)$.

We shall denote $\Sigma_R(x):=\partial\Omega\cap B_R(x)$
for any $R>0$.
We shall abbreviate $\Omega_R=\Omega_R(x)$ and
$\Sigma_R=\Sigma_R(x)$ if the point $x$ is well understood in the context.

\begin{lemma}[Boundary Poincar\'e inequality]\label{lem:G-02}
Assume that $\abs{B_R\setminus\Omega}\ge \theta\abs{B_R}$
for some $\theta >0$.
Then, for any $u\in W^{1,2}(\Omega_R)$ satisfying $u=0$ on $\Sigma_R$,
we have the following estimate:
\begin{equation}
\label{eqG-62}
\norm{u}_{L^2(\Omega_R)}\le \frac{1}{\theta} R\norm{D u}_{L^2(\Omega_R)}.
\end{equation}
\end{lemma}
\begin{proof}
Since $u=0$ in $\Sigma_R$, we may extend $u$ to a $W^{1,2}(B_R)$
function by setting $u=0$ in $S:=B_R\setminus\Omega$.
Note that $D u=0$ in $S$.
Then the lemma follows from (7.45) in \cite[p. 164]{GT}.
\end{proof}

\begin{lemma}[Boundary Caccioppoli inequality]\label{lem:G-03}
Let the operator $L$ satisfy conditions \eqref{eqP-02}, \eqref{eqP-03}.
Suppose $\vec{u}$ is a $W^{1,2}(\Omega_R)^N$ solutions of $L\vec{u}=0$
in $\Omega_R$ satisfying $\vec{u}=0$ on $\Sigma_R$.
Then, we have
\begin{equation}
\norm{D\vec{u}}_{L^2(\Omega_r)}\le \frac{C}{R-r} \norm{u}_{L^2(\Omega_R)},
\quad\forall 0<r<R,
\end{equation}
where $C=C(n,N,\lambda,\Lambda)>0$.
\end{lemma}
\begin{proof}
It is well known.
\end{proof}

\begin{definition}
We say that $\Omega$ satisfies the condition (S)
at a point $\bar{x}\in\partial\Omega$ if there exist $\theta>0$ and
$R_a\in (0,\infty]$ such that
\begin{equation}
\label{eqG-64}
\abs{B_R(\bar{x})\setminus\Omega}\ge \theta\abs{B_R(\bar{x})},\quad
\forall R< R_a.
\end{equation}
We say that $\Omega$ satisfies the condition (S) uniformly on 
$\Sigma\subset\partial\Omega$
if there exist $\theta>0$ and $R_a$ such that \eqref{eqG-64} holds
for all $\bar{x}\in\Sigma$.
\end{definition}

\begin{definition}
Let $\Omega$ satisfy the condition (S) at $\bar{x}\in\partial\Omega$.
We shall say that an operator $L$ satisfies the property (BH) if
there exist $\mu_1, H_1>0$ such that
if $\vec{u}\in W^{1,2}(\Omega_{R}(\bar{x}))^N$ is a weak solution
of the problem, $L\vec{u}=0$ in $\Omega_R(\bar{x})$ and $\vec{u}=0$
on $\Sigma_R(\bar{x})$, where $R< R_a$,
then $\vec{u}$ satisfies the following estimates:
\begin{equation}
\label{eqG-65}
\int_{\Omega_r(\bar{x})}\abs{D \vec{u}}^2 \le H_1\left(\frac{r}{s}\right)^{n-2+2\mu_1}
\int_{\Omega_s(\bar{x})}\abs{D \vec{u}}^2,\quad \forall 0<r<s\le R.
\end{equation}
\end{definition}

\begin{lemma} \label{lem:G-06}
There exists $\epsilon_0=\epsilon_0(n,\lambda_0,\Lambda_0)>0$
such that if the coefficients of the operator $L$ in \eqref{eqP-01}
satisfies \eqref{eqP-08} in Lemma~\ref{lem:P-02}, then
$L$ satisfies the property (BH)
with $\mu_1=\mu_1(n,\lambda_0,\Lambda_0,\theta)>0$
and $H_1=H_1(n,N,\lambda_0,\Lambda_0,\theta)>0$.
\end{lemma}
\begin{proof}
Throughout the proof, we shall
abbreviate $\Omega_r=\Omega_r(\bar{x})$ for any $r>0$,
$\Sigma_r=\Sigma_r(\bar{x})$, the point $\bar{x}\in\partial\Omega$ to be understood.
For any $s \le R<R_a$,
let $v^i$ ($i=1,\ldots,N$) be a unique $W^{1,2}(\Omega_s)$
solution of $L_0 v^i=0$ in $\Omega_s$ satisfying
$v^i-u^i\in W^{1,2}_0(\Omega_s)$, where
$L_0 v^i= -D_\alpha(a^{\alpha\beta}D_\beta v^i)$.

We claim that there exist $\mu_2(n,\lambda_0,\Lambda_0,\theta)>0$
and $C(n,\lambda_0,\Lambda_0,\theta)>0$ such that the following estimate holds:
\begin{equation}
\label{eqG-66}
\int_{\Omega_r}\abs{D\vec{v}}^2
\le C\left(\frac{r}{s}\right)^{n-2+2\mu_2}\int_{\Omega_s}\abs{D\vec{v}}^2,
\quad \forall 0<r<s.
\end{equation}
We first note that we may assume that $r\le s/8$;
otherwise \eqref{eqG-66} becomes trivial.
Since each $v^i$ satisfies $v^i=0$ on $\Sigma_s$,
it follows from Theorem 8.27 \cite[pp. 203--204]{GT} and Theorem 8.25
\cite[pp. 202--203]{GT} that there is
$\mu_2=\mu_2(n,\lambda_0,\Lambda_0,\theta)>0$
and $C=C(n,\lambda_0,\Lambda_0,\theta)>0$ such that
\begin{equation}
\label{eqG-67}
\osc_{\Omega_{2r}} v^i \le C r^{\mu_2}s^{-\mu_2}\sup_{\Omega_{s/4}} |v^i|
\le C r^{\mu_2} s^{-\mu_2-n/2}\Norm{v^i}_{L^2(\Omega_{s/2})}.
\end{equation}
In particular, the estimate \eqref{eqG-67} implies
$v^i(\bar{x})=\lim\limits_{x\to \bar{x}}v^i(x)=0$.
Then, Lemma~\ref{lem:G-03} and Lemma~\ref{lem:G-02} imply
that for all $i=1,\ldots,N$ (recall $r<s/8$)
\begin{equation*}
\begin{split}
\int_{\Omega_r}\abs{D v^i}^2
&\le C r^{-2}\int_{\Omega_{2r}}\Abs{v^i}^2=
C r^{-2}\int_{\Omega_{2r}}\Abs{v^i-v^i(\bar{x})}^2 \\
&\le C r^{n-2}\left(\osc_{\Omega_{2r}}\,v^i\right)^2
\le C\left(\frac{r}{s}\right)^{n-2+2\mu_2}
s^{-2}\int_{\Omega_{s/2}}\Abs{v^i}^2\\
&\le C\left(\frac{r}{s}\right)^{n-2+2\mu_2}\int_{\Omega_{s}}\Abs{D v^i}^2,
\end{split}
\end{equation*}
and thus we have proved the claim.

Next, note that $\vec{w}:=\vec{u}-\vec{v}$ belongs to $W^{1,2}_0(\Omega_s)^N$ 
and thus it satisfies
\begin{equation*}
\lambda \int_{\Omega_s} \abs{D\vec{w}}^2
\le \int_{\Omega_s} a^{\alpha\beta} D_\beta w^i D_\alpha w^i=
\int_{\Omega_s} (a^{\alpha\beta}\delta_{ij}-A^{\alpha\beta}_{ij})
D_\beta u^j D_\alpha w^i.
\end{equation*}
Therefore, we have
\begin{equation}
\label{eqG-70}
\int_{\Omega_s} \abs{D\vec{w}}^2 \le (\lambda^{-1}\norm{\epsilon}_{L^\infty})^2
\int_{\Omega_s} \abs{D\vec{u}}^2,
\end{equation}
where $\epsilon(x)$ is as defined in \eqref{eqP-08}.
By combining \eqref{eqG-66} and \eqref{eqG-70}, we obtain
\begin{equation*}
\int_{\Omega_r}\abs{D\vec{u}}^2
\le C\left(\frac{r}{s}\right)^{n-2+2\mu_2}\int_{\Omega_s}\abs{D\vec{u}}^2+
C_0 \norm{\epsilon}_{L^\infty}^2 \int_{\Omega_s} \abs{D\vec{u}}^2,
\quad \forall 0<r<s.
\end{equation*}
Now, choose a $\mu_1\in (0,\mu_2)$.
Then, from a well known iteration argument
(see, e.g., \cite[Lemma~2.1, p. 86]{Gi83}), it follows that
there is $\epsilon_0$ such that if $\norm{\epsilon}_{L^\infty}< \epsilon_0$,
then \eqref{eqG-65} holds.
\end{proof}

\begin{theorem}\label{thm:G-07}
Let the operator $L$ satisfy the properties (H) and (BH).
Assume that $\Omega$ satisfies the condition (S) at $\bar{x}\in\partial\Omega$
with parameters $\theta, R_a$.
Let $x\in\Omega$ such that
$\abs{x-\bar{x}}=d_x\le R/2$, where
$R< R_a$ is given.
Then, any weak solution  $\vec{u}$ of $L\vec{u}=0$
in $\Omega_R(\bar{x})$ satisfying $\vec{u}=0$ on $\Sigma_R(\bar{x})$, we have
\begin{equation}
\label{eqG-72}
\abs{\vec{u}(x)}\le C d_x^\mu
R^{1-n/2-\mu}\norm{D\vec{u}}_{L^2(\Omega_R(\bar{x}))},
\quad d_x:=\dist(x,\partial\Omega),
\end{equation}
where $C=C(n,N,\lambda,\Lambda,\theta,\mu_0,\mu_1,H_0,H_1)>0$ and
$\mu=\min(\mu_0,\mu_1)$.
\end{theorem}
\begin{proof}
The proof is an adaptation of a technique due to Campanato \cite{C63}.
In this proof, we shall use the notation
$\vec{u}_{x,r}:=\dashint_{\Omega_r(x)}\vec{u}$.
Also, we shall abbreviate $d=d_x$.
Observe that 
\begin{equation}
\label{eqG-73}
\Omega_d(x)=B_d(x)\subset\Omega_{2d}(x) \cap\Omega_{2d}(\bar{x}).
\end{equation}
We may assume that $R>3d$ so that $\Omega_{2d}(x)\subset \Omega_R(\bar{x})$;
otherwise $2d\le R\le 3d$ and
\eqref{eqG-72} follows from Lemma~\ref{lem:P-04}.
We estimate $\vec{u}(x)$ by
\begin{equation*}
\abs{\vec{u}(x)}
\le \abs{\vec{u}(x)-\vec{u}_{x,2d}}+
\abs{\vec{u}_{x,2d}-\vec{u}_{\bar{x},2d}}+
\abs{\vec{u}_{\bar{x},2d}} :=I+II+III.
\end{equation*}
We shall estimate $I$ first.
For any $r_1<r_2\le 2d$, we estimate
\begin{equation}
\label{eqG-75}
\abs{\vec{u}_{x,r_1}-\vec{u}_{x,r_2}}^2\le
2 \abs{\vec{u}(z)-\vec{u}_{x,r_1}}^2+2\abs{\vec{u}(z)-\vec{u}_{x,r_2}}^2.
\end{equation}
Note that since $B_d(x)\subset \Omega$, we have
\begin{equation*}
\abs{\Omega_r(x)}\ge C r^n,\quad \forall r\le 2d.
\end{equation*}
Therefore, by integrating \eqref{eqG-75} over $\Omega_{r_1}(x)$
with respect to $z$, we estimates
\begin{equation}
\label{eqG-77}
\abs{\vec{u}_{x,r_1}-\vec{u}_{x,r_2}}^2\le
C r_1^{-n} \left(\int_{\Omega_{r_1}} \abs{\vec{u}-\vec{u}_{x,r_1}}^2
+\int_{\Omega_{r_2}}\abs{\vec{u}-\vec{u}_{x,r_2}}^2\right).
\end{equation}
Since $\vec{u}=0$ on $\Sigma_R(\bar{x})$, we may extend
$u$ to $B_R(\bar{x})$ as a $W^{1,2}$ function by setting
$\vec{u}=0$ on $B_R(\bar{x})\setminus\Omega$.
Therefore, by a version of Poincar\'e inequality
(see, e.g. (7.45) in \cite[p. 164]{GT}), we have for all $r\le 2d$,
\begin{equation}
\label{eqG-78}
\int_{\Omega_{r}} \abs{\vec{u}-\vec{u}_{x,r}}^2
\le \int_{B_{r}} \abs{\vec{u}-\vec{u}_{x,r}}^2 \le
C r^2 \int_{B_{r}} \abs{D\vec{u}}^2
=C r^2 \int_{\Omega_{r}} \abs{D\vec{u}}^2.
\end{equation}
Therefore, by \eqref{eqG-77} and \eqref{eqG-78}, we obtain
\begin{equation}
\label{eqG-79}
\abs{\vec{u}_{x,r_1}-\vec{u}_{x,r_2}}^2\le
C r_1^{-n} \left(r_1^2\int_{\Omega_{r_1}(x)} \abs{D\vec{u}}^2
+r_2^2\int_{\Omega_{r_2}(x)}\abs{D\vec{u}}^2\right).
\end{equation}
Next, we claim that the following estimate holds:
\begin{equation}
\label{eqG-80}
\int_{\Omega_r(x)}\abs{D \vec{u}}^2 \le
C\left(\frac{r}{R}\right)^{n-2+2\mu}
\int_{\Omega_{R(\bar{x})}}\abs{D\vec{u}}^2, \quad \forall r\le 2d.
\end{equation}
We first consider the case when $r\le d$.  Note that in this case,
we have $\Omega_r(x)= B_r(x)$ and $\Omega_d(x)= B_d(x)$.
Since $L$ satisfies (H), it follows from \eqref{eqG-73} that
\begin{equation}
\label{eqG-81}
\int_{\Omega_r(x)}\abs{D \vec{u}}^2
\le C\left(\frac{r}{d}\right)^{n-2+2\mu}
\int_{\Omega_d(x)}\abs{D \vec{u}}^2
\le C\left(\frac{r}{d}\right)^{n-2+2\mu}
\int_{\Omega_{2d}(\bar{x})}\abs{D \vec{u}}^2.
\end{equation}
On the other hand, since $L$ satisfies (BH), it follows
from \eqref{eqG-65} that
\begin{equation}
\label{eqG-82}
\int_{\Omega_{2d}(\bar{x})}\abs{D \vec{u}}^2
\le C\left(\frac{d}{R}\right)^{n-2+2\mu}
\int_{\Omega_{R}(\bar{x})}\abs{D \vec{u}}^2.
\end{equation}
By combining \eqref{eqG-81} and \eqref{eqG-82}, we obtain \eqref{eqG-80}.
Next, consider the case when $d<r$.
In this case, we have $\Omega_r(x)\subset \Omega_{2r}(\bar{x})$,
and thus it follows from \eqref{eqG-65}
\begin{equation*}
\int_{\Omega_r(x)}\abs{D \vec{u}}^2
\le \int_{\Omega_{2r}(\bar{x})}\abs{D \vec{u}}^2
\le C\left(\frac{r}{R}\right)^{n-2+2\mu}
\int_{\Omega_{R}(\bar{x})}\abs{D \vec{u}}^2.
\end{equation*}
We proved the claim \eqref{eqG-80}.

Now, by using \eqref{eqG-80}, we estimates \eqref{eqG-79} as follows
(recall $r_1<r_2\le 2d$):
\begin{equation}
\label{eqG-84}
\abs{\vec{u}_{x,r_1}-\vec{u}_{x,r_2}}^2\le
C r_1^{-n}(r_1^{n+2\mu}+r_2^{n+2\mu})
R^{2-n-2\mu} \int_{\Omega_{R}(\bar{x})}\abs{D\vec{u}}^2.
\end{equation}
For any $r\le 2d$, set $r_1=r 2^{-(i+1)}$ and $r_2=r 2^{-i}$
in \eqref{eqG-84} to get
\begin{equation*}
\abs{\vec{u}_{x,r2^{-(i+1)}}-\vec{u}_{x,r2^{-i}}}^2
\le C r^{2\mu} 2^{-2\mu(i+1)}
R^{2-n-2\mu} \int_{\Omega_{R}(\bar{x})}\abs{D \vec{u}}^2.
\end{equation*}
Therefore, for $0\le j<k$, we obtain
\begin{equation}
\label{eqG-86}
\begin{split}
\abs{\vec{u}_{x,r2^{-k}}-\vec{u}_{x,r2^{-j}}}
&\le \sum_{i=j}^{k-1} \abs{\vec{u}_{x,r2^{-(i+1)}}-\vec{u}_{x,r2^{-i}}}\\
&\le C r^{\mu} \left(\sum_{i=j}^\infty 2^{-\mu(i+1)}\right) R^{1-n/2-\mu}
\norm{D \vec{u}}_{L^2(\Omega_{R}(\bar{x}))}\\
&=C 2^{-j\mu} r^{\mu} R^{1-n/2-\mu}
\norm{D \vec{u}}_{L^2(\Omega_{R}(\bar{x}))}.
\end{split}
\end{equation}
By setting $r=2d$, $j=0$, and letting $k\to\infty$ in \eqref{eqG-86},
we obtain
\begin{equation}
\label{eqG-87}
I=\abs{\vec{u}(x)-\vec{u}_{x,2d}}
\le C d^{\mu} R^{1-n/2-\mu}
\norm{D \vec{u}}_{L^2(\Omega_{R}(\bar{x}))}.
\end{equation}

Next, we estimate $III$.
Since $\abs{B_r(\bar{x})\cap B_d(x)}\ge C r^n$ for $r\le 2d$, we have
\begin{equation}
\label{eqG-88}
\abs{\Omega_r(\bar{x})}\ge C r^n,\quad \forall r\le 2d.
\end{equation}
Also, as in \eqref{eqG-78}, we have for all $r\le 2d$
(recall $\vec{u}\equiv 0$ on $B_R(\bar{x})\setminus\Omega$)
\begin{equation}
\label{eqG-89}
\int_{\Omega_{r}} \abs{\vec{u}-\vec{u}_{\bar{x},r}}^2
\le \int_{B_{r}} \abs{\vec{u}-\vec{u}_{\bar{x},r}}^2 \le
C r^2 \int_{B_{r}} \abs{D\vec{u}}^2
=C r^2 \int_{\Omega_{r}} \abs{D\vec{u}}^2.
\end{equation}
Therefore, as in \eqref{eqG-79} we have
for $r_1<r_2\le 2d$,
\begin{equation*}
\abs{\vec{u}_{\bar{x},r_1}-\vec{u}_{\bar{x},r_2}}^2\le
C r_1^{-n} \left(r_1^2\int_{\Omega_{r_1}(\bar{x})} \abs{D\vec{u}}^2
+r_2^2\int_{\Omega_{r_2}(\bar{x})}\abs{D\vec{u}}^2\right).
\end{equation*}
Then, by using the property (BH), we obtain
(c.f. \eqref{eqG-86}, \eqref{eqG-87})
\begin{equation}
\label{eqG-91}
\abs{\hat{\vec{u}}(\bar{x})-\vec{u}_{\bar{x},2d}}
\le C d^{\mu} R^{1-n/2-\mu}
\norm{D \vec{u}}_{L^2(\Omega_{R}(\bar{x}))},
\end{equation}
where $\hat{\vec{u}}(\bar{x}):=\lim_{k\to\infty} \vec{u}_{\bar{x},2^{-k}r}$.
(note that \eqref{eqG-86} implies $\hat{\vec{u}}(\bar{x})$ exists).
It follows from \eqref{eqG-88}, \eqref{eqG-62}, and \eqref{eqG-65}
that for any $r\le 2d$, we have
\begin{equation*}
\begin{split}
\abs{\vec{u}_{\bar{x},r}}^2
&\le \dashint_{\Omega_r(\bar{x})}\abs{\vec{u}}^2\le
Cr^{-n}\int_{\Omega_r(\bar{x})}\abs{\vec{u}}^2 \\
&\le Cr^{2-n}\int_{\Omega_r(\bar{x})}\abs{D\vec{u}}^2
\le C r^{2-n}\left(\frac{r}{R}\right)^{n-2+2\mu}\int_{\Omega_R(\bar{x})}
\abs{D \vec{u}}^2\\
&= C r^{2\mu}R^{2-n-2\mu}\int_{\Omega_R(\bar{x})} \abs{D \vec{u}}^2,
\end{split}
\end{equation*}
and thus that $\hat{\vec{u}}(\bar{x})=0$.
Therefore, by \eqref{eqG-91} we obtain
\begin{equation}
\label{eqG-93}
III= \abs{\vec{u}_{\bar{x},2d}}
=\abs{\hat{\vec{u}}(\bar{x})-\vec{u}_{\bar{x},2d}}
\le C d^{\mu} R^{1-n/2-\mu}
\norm{D \vec{u}}_{L^2(\Omega_{R}(\bar{x}))}.
\end{equation}

Finally, we estimate $II$.
\begin{equation}
\label{eqG-94}
\abs{\vec{u}_{x,2d}-\vec{u}_{\bar{x},2d}}^2 \le
2\abs{\vec{u}(z)-\vec{u}_{x,2d}}^2+ 2\abs{\vec{u}(z)-\vec{u}_{\bar{x},2d}}^2.
\end{equation}
By integrating \eqref{eqG-94} over
$B_d(x)\subset \Omega_{2d}(x)\cap \Omega_{2d}(\bar{x})$
with respect to $z$, we estimate
\begin{equation}
\label{eqG-95}
\begin{split}
\abs{\vec{u}_{x,2d}-\vec{u}_{\bar{x},2d}}^2 &\le
C d^{-n}\left(\int_{\Omega_{2d}(x)}\abs{\vec{u}-\vec{u}_{x,2d}}^2+
\int_{\Omega_{2d}(\bar{x})}\abs{\vec{u}-\vec{u}_{\bar{x},2d}}^2\right)\\
&\le C d^{2-n}\left(\int_{\Omega_{2d}(x)}\abs{D\vec{u}}^2+
\int_{\Omega_{2d}(\bar{x})}\abs{D\vec{u}}^2\right)\\
&\le C d^{2\mu}R^{2-n-2\mu}\int_{\Omega_{R}(\bar{x})}\abs{D\vec{u}}^2,
\end{split}
\end{equation}
where we have used \eqref{eqG-78}, \eqref{eqG-89}, \eqref{eqG-80},
and \eqref{eqG-65}.
Therefore, by combining \eqref{eqG-87}, \eqref{eqG-93},
and \eqref{eqG-95}, we obtain \eqref{eqG-72}.
\end{proof}

\begin{theorem}\label{thm:G-08}
Let the operators $L$, ${}^t\!L$ satisfy the properties (H) and (BH).
Assume that $\Omega$ satisfies the condition (S) uniformly on $\partial\Omega$
with parameters $\theta, R_a$.
Denote
\begin{equation*}
R_{x,y}:=\min(\abs{x-y}, 4R_a).
\end{equation*}
Then the Green matrix $\vec{G}(x,y)$ satisfies
\begin{gather}
\label{eqG-97}
\abs{\vec{G}(x,y)}\le C d_x^\mu R_{x,y}^{1-n/2-\mu} d_y^{1-n/2}
\quad\text{if } d_x\le R_{x,y}/8,\\
\label{eqG-98}
\abs{\vec{G}(x,y)}\le C d_y^\mu R_{x,y}^{1-n/2-\mu} d_x^{1-n/2}
\quad\text{if } d_y\le R_{x,y}/8,
\end{gather}
where $C=C(n,N,\lambda,\Lambda,\theta,\mu_0,\mu_1,H_0,H_1)>0$ and
$\mu=\min(\mu_0,\mu_1)$.
As a consequence, we have $\vec{G}(\,\cdot\,,y)=0$,
$\vec{G}(x,\,\cdot\,)=0$ on $\partial\Omega$ in the usual sense.
\end{theorem}
\begin{proof}
We only need to prove \eqref{eqG-97}, for \eqref{eqG-98} will then follow
from \eqref{eqG-45}.
Set $R=R_{x,y}/4$, $r=d_y/2$,
and choose $\bar{x}\in\partial\Omega$ such that $\abs{x-\bar{x}}=d_x$.
Then, since
\begin{equation*}
d_y\le \abs{x-y}+d_x\le\tfrac{9}{8}\abs{x-y},
\end{equation*}
we have
\begin{equation*}
\abs{y-\bar{x}}\ge \abs{x-y}-d_x\ge \tfrac{7}{8}\abs{x-y}\ge R+r,
\end{equation*}
and thus,
$\Omega_R(\bar{x})\subset \Omega\setminus B_r(y)$.
Now, we apply Theorem~\ref{thm:G-07} with $\vec{u}=\vec{G}(\,\cdot\,,y)$.
Then, by \eqref{eqG-72} and \eqref{eqG-34},
we obtain
\begin{equation*}
\abs{\vec{G}(x,y)}
\le C d_x^\mu R^{1-n/2-\mu}
\norm{D\vec{G}(\,\cdot\,,y)}_{L^2(\Omega\setminus B_r(y))}
\le C d_x^\mu R_{x,y}^{1-n/2-\mu} d_y^{1-n/2}.
\end{equation*}
The proof is complete.
\end{proof}

\begin{remark}
We note that in the scalar case, the maximum principle yields
(see \cite[Theorem 1.1]{GW})
\begin{equation}
\label{eqG-104}
G(x,y)\le C \abs{x-y}^{2-n},\quad \forall x\neq y\in\Omega.
\end{equation}
Then, by the boundary Caccioppoli inequality, we have
(c.f. \eqref{eqG-11}--\eqref{eqG-14})
\begin{equation*}
\int_{\Omega\setminus B_r(y)} \abs{D G(\,\cdot\,,y)}^2 \le C r^{2-n},
\quad \forall r>0.
\end{equation*}
Therefore, in the scalar case we don't need to require that $r< d_y/2$
(or $r< d_x/2$) in the proof of Theorem~\ref{thm:G-08} and we may as well
set $r= \abs{x-y}/2$ to get
\begin{gather*}
G(x,y)\le C d_x^\mu R_{x,y}^{1-n/2-\mu} \abs{x-y}^{1-n/2}
\quad\text{if }d_x\le R_{x,y}/8,\\
G(x,y)\le C d_y^\mu R_{x,y}^{1-n/2-\mu} \abs{x-y}^{1-n/2}
\quad\text{if }d_y\le R_{x,y}/8.
\end{gather*}
In particular, if $G(x,y)$ is the Green's function on $\bR^n_{+}$, then
we obtain
\begin{gather*}
G(x,y)\le C d_x^\mu \abs{x-y}^{2-n-\mu} \quad\text{if }d_x\le \abs{x-y}/8,\\
G(x,y)\le C d_y^\mu \abs{x-y}^{2-n-\mu} \quad\text{if }d_y\le \abs{x-y}/8,
\end{gather*}
for $\partial\bR^n_{+}$ satisfies the condition (S) with $\theta=1/2$ and
$R_a=\infty$.
\end{remark}

\section{Remarks on VMO coefficients case}\label{sec:V}
\begin{definition}[Sarason \cite{Sarason}]\label{def:V-01}
For a measurable function $f$ defined on $\bR^n$,
we shall denote $\overline{f}_{x,r}=\dashint_{B_r(x)}f$
and for $0<\delta<\infty$ we define
\begin{equation}
\label{eqV-01}
M_\delta(f):=\sup_{x\in\bR^n}\,\sup_{r\le \delta}\,\dashint_{B_r(x)}
\abs{f-\overline{f}_{x,r}};
\quad M_0(f):=\lim_{\delta\to 0} M_\delta(f).
\end{equation}
We shall say that $f$ belongs to $\VMO$ if
$M_0(f)=0$.
\end{definition}

\begin{definition}\label{def:V-02}
We say that the operator $L$ satisfies the property $\Hloc$ if
there exist $\mu_0,H_0,R_c>0$ such that all weak solutions $\vec{u}$ of
$L\vec{u}=0$ in $B_R=B_R(x_0)$ with $R\le R_c$ satisfy
\begin{equation}
\label{eqV-02}
\int_{B_r}\abs{D \vec{u}}^2 \le H_0 \left(\frac{r}{s}\right)^{n-2+2\mu_0}
\int_{B_s}\abs{D \vec{u}}^2, \quad 0<r<s\le R.
\end{equation}
Similarly, we say that the transpose operator ${}^t\!L$
satisfies the property $\Hloc$ if corresponding
estimates hold for all weak solutions $\vec{u}$ of
${}^t\!L\vec{u}=0$ in $B_R$ with $R\le R_c$.
\end{definition}

\begin{lemma}\label{lem:V-03}
Let the coefficients of the operator $L$ in \eqref{eqP-01}
satisfy the conditions \eqref{eqP-02} and \eqref{eqP-03}.
If the coefficients belong to $\VMO$ in addition, then 
$L$ satisfies the property $\Hloc$.
\end{lemma}
\begin{proof}
It is well known that if the coefficients are uniformly continuous, then
$L$ satisfies the property $\Hloc$; see e.g. \cite[pp. 87--89]{Gi83}.
Essentially, the same proof carries over to the $\VMO$ coefficients case.
One only needs to make a note of the following two facts.
First, a theorem of Meyers \cite{Meyers} implies that
there is some $p=p(n,N,\lambda,\Lambda)>2$ such that
if $\vec{u}$ is a weak solution of $L\vec{u}=0$ in $B_R(x)$, then
\begin{equation*}
\left(\dashint_{B_r(x)}\abs{D\vec{u}}^p\right)^{1/p}
\le C\left(\dashint_{B_{2r(x)}}\abs{D\vec{u}}^2\right)^{1/2},
\quad \forall r<R/2.
\end{equation*}
Secondly, note that the John-Nirenberg theorem \cite{JN} implies that
\begin{equation*}
\left(\dashint_{B_r(x)}\abs{f-\overline{f}_{r,x}}^q\right)^{1/q}
\le C(n,q) M_\delta(f),\quad \forall r<c(n)\delta,\quad \forall q\in (0,\infty),
\end{equation*}
where $M_\delta(f)$ is defined as in \eqref{eqV-01}.
For the details, we refer to \cite[pp. 47--48]{AT}.
\end{proof}

In the rest of this section, we shall assume that the operators $L$ and
${}^t\!L$ satisfy the property $\Hloc$ with parameters $\mu_0, H_0, R_c$.
We shall denote
\begin{equation}
\label{eqV-05}
r_x:=\min(d_x,R_c), \quad \bar{r}_{x,y}:=\min(\bar{d}_{x,y},R_c),
\end{equation}
where $d_x=\dist(x,\partial\Omega)$ and $\bar{d}_{x,y}$ is as in
\eqref{eqG-39}.
It is routine to check that all estimates appearing in Section~\ref{sec:G-01}
remain valid if $d_x$, $\bar{d}_{x,y}$ are replaced by $r_x$,
$\bar{r}_{x,y}$, respectively.
Therefore, we have the following theorem:

\begin{theorem}
Let $\Omega$ be an open connected set in $\bR^n$. Denote
$d_x:=\dist(x,\partial\Omega)$ for $x\in\Omega$; we set $d_x=\infty$ if 
$\Omega=\bR^n$.
Assume that operators $L$ and ${}^t\!L$ satisfy the property $\Hloc$.
Then, there exists a unique Green's matrix
$\vec{G}(x,y)=(G_{ij}(x,y))_{i,j=1}^N$ ($x,y\in\Omega, x\neq y$)
which is continuous in $\set{(x,y)\in\Omega\times\Omega:x\neq y}$
and such that $\vec G(x,\,\cdot\,)$ is locally integrable in
$\Omega$ for all $x\in\Omega$ and that for all
$\vec{f}=(f^1,\ldots,f^N)^T \in C^\infty_c(\Omega)^N$,
the function $\vec{u}=(u^1,\ldots,u^N)^T$ given by
\begin{equation}
\label{eqV-70}
\vec{u}(x):=\int_{\Omega} \vec{G}(x,y)\vec{f}(y)\,dy
\end{equation}
belongs to $Y^{1,2}_0(\Omega)^N$ and satisfies $L\vec{u}=\vec{f}$ in the sense
\begin{equation}
\label{eqV-71}
\int_{\Omega} A^{\alpha\beta}_{ij} D_\beta u^j D_\alpha \phi^i = \int_{\Omega}
f^i \phi^i,\quad \forall \vec{\phi}\in C^\infty_c(\Omega)^N.
\end{equation}
Moreover, $\vec{G}(x,y)$ has the properties that
\begin{equation}
\label{eqV-54}
\int_{\Omega}A^{\alpha\beta}_{ij} D_\beta G_{jk}(\,\cdot\,,y)
D_\alpha \phi^i = \phi^k(y),
\quad \forall \vec{\phi}\in C^\infty_c(\Omega)^N
\end{equation}
and that for all $\eta\in C^\infty_c(\Omega)$ satisfying
$\eta\equiv 1$ on $B_r(y)$ for some $r<d_y$,
\begin{equation}
(1-\eta)\vec{G}(\,\cdot\,,y)\in Y^{1,2}_0(\Omega)^{N\times N}.
\end{equation}
Furthermore, $\vec{G}(x,y)$ satisfies the following estimates:
For $r_x, r_y, \bar{r}_{x,y}$ as in \eqref{eqV-05},
\begin{gather}
\norm{\vec{G}(\,\cdot\,,y)}_{L^p(B_r(y))}\le C_p\,r^{2-n+n/p},
\quad \forall r< r_y,
\quad \forall p\in [1,\tfrac{n}{n-2}),\\
\norm{\vec{G}(x,\,\cdot\,)}_{L^p(B_r(x))}\le C_p\,r^{2-n+n/p},
\quad \forall r< r_x,
\quad \forall p\in [1,\tfrac{n}{n-2}),
\end{gather}
\begin{gather}
\norm{D\vec{G}(\,\cdot\,,y)}_{L^p(B_r(y))}\le C_p\,r^{1-n+n/p},
\quad \forall r< r_y,
\quad \forall p\in [1,\tfrac{n}{n-1}),\\
\norm{D\vec{G}(x,\,\cdot\,)}_{L^p(B_r(x))}\le C_p\,r^{1-n+n/p},
\quad \forall r< r_x,
\quad \forall p\in [1,\tfrac{n}{n-1}),
\end{gather}
\begin{gather}
\norm{\vec{G}(\,\cdot\,,y)}_{Y^{1,2}(\Omega\setminus B_r(y))}
\le C r^{1-n/2},\quad \forall r< r_y/2,\\
\norm{\vec{G}(x,\,\cdot\,)}_{Y^{1,2}(\Omega\setminus B_r(x))}
\le C r^{1-n/2},\quad \forall r< r_x/2,
\end{gather}
\begin{gather}
\abs{\set{x\in\Omega:\abs{\vec{G}(x,y)}>t}}\le C t^{-\frac{n}{n-2}},
\quad\forall t> (r_y/2)^{2-n},\\
\abs{\set{y\in\Omega:\abs{\vec{G}(x,y)}>t}}\le C t^{-\frac{n}{n-2}},
\quad\forall t> (r_x/2)^{2-n},
\end{gather}
\begin{gather}
\abs{\set{x\in\Omega:\abs{D_x \vec{G}(x,y)}>t}}\le C t^{-\frac{n}{n-1}},
\quad\forall t> (r_y/2)^{1-n},\\
\abs{\set{y\in\Omega:\abs{D_y \vec{G}(x,y)}>t}}\le C t^{-\frac{n}{n-1}},
\quad\forall t> (r_x/2)^{1-n},
\end{gather}
\begin{equation}
\label{eqV-18}
\abs{\vec{G}(x,y)}\le C \bar{r}_{x,y}^{2-n},\quad
\forall x,y\in \Omega,
\end{equation}
\begin{gather}
\label{eqV-19}
\abs{\vec{G}(x,y)-\vec{G}(z,y)}\le
C \abs{x-z}^{\mu_0} \bar{r}_{x,y}^{2-n-\mu_0}\quad\text{if }
\abs{x-z}<\bar{r}_{x,y}/2,\\
\abs{\vec{G}(x,y)-\vec{G}(x,z)}\le
C \abs{y-z}^{\mu_0} \bar{r}_{x,y}^{2-n-\mu_0}\quad\text{if }
\abs{y-z}<\bar{r}_{x,y}/2,
\end{gather}
where $C=C(n,N,\lambda,\Lambda,\mu_0,H_0)>0$ and 
$C_p=C_p(n,N,\lambda,\Lambda,\mu_0,H_0,p)>0$.
\end{theorem}

\begin{remark}
Dolzmann-M\"{u}ller \cite{DM} derived a global estimate
\begin{equation}
\label{eqV-21}
\abs{\vec{G}(x,y)}\le C \abs{x-y}^{2-n}\quad \forall x,y\in \Omega,
\quad x\neq y,
\end{equation}
assuming that $\Omega$ is a bounded $C^1$ domain.
We have not attempted to derive the corresponding estimate here.
However, we would like to point out that the constant $C$ in their estimate depends on the domain (e.g., the diameter of the
domain and also some characteristics of $\partial\Omega$) while our interior estimate \eqref{eqV-18} does not.
\end{remark}

\end{document}